\documentclass[twoside,11pt]{article}

\usepackage{amsfonts}
\usepackage{mathrsfs}
\usepackage{amsmath,amsfonts}

\newcommand{\D}{\displaystyle}

\begin{document}

\title{{\bf  Convergence Rates to Stationary Solutions of a Gas-liquid Model with External Forces and Vacuum}}
\author{{L{\sc ong} F{\sc an}\ \ \ \ \ \ \ Q{\sc ingqing} L{\sc iu} \ \ \ \ \ \ C{\sc hangjiang} Z{\sc hu}\thanks{Corresponding author.
Email: cjzhu@mail.ccnu.edu.cn}}   \\  \\
{The Hubei Key Laboratory of Mathematical Physics} \\
{School of Mathematics and Statistics }
\\{Central China Normal University, Wuhan 430079, P.R. China}}
\date{}
\maketitle

 \textbf{{\bf Abstract:}} In this paper, we study the
asymptotic behavior of solutions to a Gas-liquid model with external
forces and general pressure law. Under some suitable assumptions on
the initial date and $\gamma>1$, if
$\theta\in(0,\frac{\gamma}{2}]\cap(0,\gamma-1]\cap(0,1-\alpha\gamma]$,
we prove the weak solution $(cQ(x,t),u(x,t))$ behavior
asymptotically to the stationary one by adapting and modifying the
technique of weighted estimates. In addition, if
$\theta\in(0,\frac{\gamma}{2}]\cap(0,\gamma-1)\cap(0,1-\alpha\gamma]$,
following the same idea in \cite{Fang-Zhang4}, we estimate the
stabilization rate of the solution as time tends to infinity in the
sense of $L^\infty$ norm.

 \bigbreak

\bigbreak \textbf{{\bf Mathematics subject classifications (2000):}}
76T10, 35L65, 35B40.

 \vspace{3mm}

\section*{Contents}

1. Introduction \dotfill 1

\noindent 2. Formulation of  problem and main results \dotfill 4

\noindent 3. Uniform {\it a priori} estimates \dotfill 6

\noindent 4. Asymptotic behavior \dotfill 19

\noindent 5. Stabilization rate estimates \dotfill 21

\noindent Acknowledgements \dotfill 28

\noindent References \dotfill 28

\vspace{4mm}

\bigbreak
\section{Introduction}

The one-dimensional two-phase model of the drift-flux type is
frequently used to simulate unsteady, compressible flow of liquid
and gas in pipes and wells (cf. \cite{Delhaye,Evje-Fjelde}). The
model consists of two
 mass conservation equations corresponding to each of the two phases gas ($g$) and liquid ($l$) and
 one equation for the momentum of the mixture and takes the following
 form
\begin{eqnarray}\label{1.1}
\left \{
\begin{array}{l}
\partial_{t}[\alpha_{g} \rho_{g}]+  \partial_{x}[\alpha_{g}
\rho_{g}u_{g}] = 0 ,\\[3mm]
\partial_{t}[\alpha_{l} \rho_{l}] +  \partial_{x}[\alpha_{l}
\rho_{l}u_{l}] =0,\\[3mm]
 \partial_{t}[\alpha_{g}\rho_{g}u_{g}+  \alpha_{l}\rho_{l}u_{l}]+\partial_{x}[\alpha_{g}\rho_{g}u_{g}^{2}+
 \alpha_{l}\rho_{l}u_{l}^{2}+p]=-q+\partial_{x}[\varepsilon
 \partial_{x}u_{mix}],
\end{array}
\right.
\end{eqnarray}
where $u_{mix}=\alpha_{g}u_{g}+\alpha_{l}u_{l}$ and the unknown
 variables $\alpha_{g},\alpha_{l} \in [0,1]$  denote volume
 fractions satisfying the fundamental relation:
\begin{eqnarray}
 \alpha_{g}+\alpha_{l}=1.
\end{eqnarray}
 Furthermore, the other unknown
 variables $\rho_{g},\ \rho_{l},\ u_{g},\ u_{l}$ denote gas density, liquid
 density, velocities of gas and liquid respectively, whereas $p$ is  the common
 pressure for both phases, $q$ presents external forces, like
 gravity and friction, and $\varepsilon>0$ denotes viscosity.

As in \cite{FE}, we assume  the liquid is incompressible,
 and the gas is polytropic,  i.e.,
\begin{eqnarray*}
\rho_{l}={\it constant},\ \ \ p=A\rho_{g}^{\gamma}, \ \ \ \gamma>1,
\ \ \ A>0.
\end{eqnarray*}

For the relation between the drift-flux model and the more general
two-fluid and the simplification of the model, one can refer to
\cite{ST} and \cite{FE} for more details. Here, we have skipped them
and  we directly study the following gas-liquid model described in
terms of Lagrangian variables (cf. \cite{FE}) with frictional force
term $-fm^2u|u|\ (f\geq 0)$ and gravity term $g (>0)$:
\begin{eqnarray}\label{equation}
\left
\{\begin{array}{l} n_t+(nm)u_x=0\\[3mm]
m_t+m^2u_x=0, \\[3mm]
      u_t+p(n,m)_x=-fm^2u|u|+g+(\varepsilon(m,n)mu_x)_x, \ \  0<x<1, \ \  t>0,\\
  \end{array}
\right.
\end{eqnarray}
with
\begin{equation}\label{PE}
p(n,m)=A\left(\frac{n}{\rho_l-m}\right)^{\gamma},\ \
\varepsilon(m,n)=\frac{B n^\theta}{(\rho_l-m)^{\theta+1}},\ \ A,\
B>0.
\end{equation}
Boundary conditions are given by
\begin{eqnarray}\label{Bound}
n(0,t)=m(0,t)=0,\ \ u(1,t)=0,\ \ t\geq0,
\end{eqnarray}
whereas initial data are
\begin{eqnarray}\label{initial}
n(x,0)=n_0(x),\ m(x,0)=m_0(x),\ u(x,0)=u_0(x),\ x \in [0,1].
\end{eqnarray}
Without loss of generality, we take $A=B=1$.

 Let's  first  review some of the previous works in this
direction. For the simplified model obtained by neglecting
frictional force and gravity in $\eqref{equation}$, there has been a
large number of research results when the viscosity function
$\varepsilon(m,n)$ takes different forms. Clearly speaking, Evje and
Karlsen in \cite{Evje-Karlsen} studied the existence of the global
weak solutions when the initial masses connected to vacuum
discontinuously for the viscosity function $\varepsilon(m,n)$ taking
the following form
\begin{equation}\label{varepsilon1}
\varepsilon(m)=\frac{m^\theta}{(\rho_l-m)^{\theta+1}}, \ \ \theta\in
(0,\frac{1}{3}).
\end{equation}
This result was later improved to the case $\theta\in (0,1]$ by Yao
and Zhu in \cite{Yao-Zhu2}. When the initial masses connected to
vacuum continuously, Evje, Flatten and Friis in
\cite{Evje-Flatten-Friis} also obtained the global existence and
uniqueness of weak solutions for the viscosity $\varepsilon(m,n)$
taking the following form
\begin{equation}\label{varepsilon2}
\varepsilon(n,m)=\frac{n^\theta}{(\rho_l-m)^{\theta+1}}, \ \
\theta\in (0,\frac{1}{3}).
\end{equation}
When the viscosity is constant and the initial masses connected to
vacuum continuously, the existence of  global weak solutions was
obtained by Yao and Zhu in \cite{Yao-Zhu1}. Recently, Friis and Evje
in \cite{FE} proved the existence of global weak solutions to the
initial boundary value problem of \eqref{equation} with initial data
\eqref{initial} and boundary conditions
$[p(n,m)-\varepsilon(m,n)mu_x](0, t) = 0,\ u(1,t)=0$ when the
viscosity $\varepsilon(m,m)$ takes the form \eqref{varepsilon1}.

However, there are few results on the asymptotic behavior and decay
rate estimates for the initial boundary value problem
\eqref{equation}-\eqref{initial}. In this paper, we rewrite our
problem into (\ref{E2})-(\ref{bc}) similar to the model of 1D
compressible Navier-Stokes equations with gravity by using the
variable transformations as in
\cite{Evje-Karlsen,Evje-Flatten-Friis} and study the asymptotic
behavior and decay rate of (\ref{E2})-(\ref{bc}) by similar line as
in \cite{Duan,Fang-Zhang4,Zi-Zhu}. According to the method in
\cite{Zi-Zhu}, we know that the uniform upper and lower bounds of
$cQ$ in \eqref{E2} and the uniform bound of
\begin{equation}\label{uxt}
\int^t_0\int^1_0c^{\theta}Q^{1+\theta}u_{xt}^2dxds
\end{equation}
play a very important role in studying the asymptotic behavior of
$cQ$ and $u$. In order to obtain these uniform bounds, we have to
deal with the difficulties which come from the additional external
force $-h(Q)u|u|$ in \eqref{E2} which is more complex than the only
constant external force $g$ in \cite{Zi-Zhu}. To overcome this
difficulty, we use \textit{a priori} assumption in the proof of the
uniform upper and lower bounds of $cQ$. Furthermore, we  obtain the
higher regularity of the velocity function $u$ by using a new skill.
Clearly speaking, we get the uniform estimate of
$\|u(\cdot,t)\|_{L^p([0,1])}$ ($ p=2,4,5$) and
$\|u\|_{L^q([0,1]\times[0,\infty))}$ ($ q=2,3,4$) with the help of
the recurrence method. The higher regularity of the velocity
function $u$ improves the estimate of
$\|u\|_{L^2([0,1]\times[0,\infty))}$ in the previous works
\cite{Duan,Fang-Zhang4,Zi-Zhu}.
 Based on these estimates, we get the uniform bound of \eqref{uxt}.

It is necessary for us  to illustrate that the main methods used to
obtain our results are similar to those  in
\cite{Duan,Fang-Zhang4,Zi-Zhu}. In view of this, let's review some
of the relevant works about Navier-Stokes equations with
density-dependent viscosity and vacuum.  For the case without
external force, Guo and Zhu in \cite{Guo1,Guo2} gave the asymptotic
behavior and decay rate of the density function $\rho(x,t)$ when the
initial density connects to vacuum continuously. Zhu in \cite{Zhu}
investigated the asymptotic behavior and decay rate estimates on the
density function $\rho(x,t)$ by overcoming some new difficulties
which came from the appearance of boundary layers when the initial
density connects to vacuum discontinuously. In \cite{Guo1,Guo2,Zhu},
the  auxiliary function $w(x,t)$ introduced by Nagasawa in
\cite{Nagasawa} was used to investigate the decay rate of
$\rho(x,t)$.  For the other case with gravity, under some
assumptions on the initial data, Zhang and Fang in
\cite{Fang-Zhang4} proved that the solution converges to the
stationary states as time goes to infinity provided
$\theta\in(0,\gamma-1)\cap(0,\frac{\gamma}{2}]$ and $\gamma>1$. The
stabilization rates were  also estimated in several norms. Duan in
\cite{Duan} generalized part result in \cite{Fang-Zhang4}, and
showed that the solution converges to the stationary state in the
sense of integral when $\gamma=2$, $\theta=1$. Recently, Zhu and Zi
in \cite{Zi-Zhu} improved the results in \cite{Duan,Fang-Zhang4} in
the sense that $\theta\in(0,\gamma-1]\cap(0,\frac{\gamma}{2}]$.

The rest of this paper is organized as follows. In Section 2, the
system is translated into a more simple one, then we give the
definition of the weak solution and state the main results. In
Section 3, we derive some crucial uniform estimates for studying the
asymptotic behavior and the decay rate estimates. In Section 4, the
asymptotic behavior of weak solution will be given. In Section 5, we
will establish stabilization rate estimates of the solution as time
tends to infinity.

\section{Formulation of problem and main results}
\setcounter{equation}{0}

In this section, we rewrite our problem
\eqref{equation}-\eqref{initial} into \eqref{E2}-\eqref{bc} by using
the same transformations as in
\cite{Evje-Karlsen,Evje-Flatten-Friis}, then we state our main
results in Theorem 2.3 and Theorem 2.4.

Firstly we introduce the transformations
\begin{eqnarray}\label{tran}
c=\frac{n}{m},\ \ \ \ Q(m)=\frac{m}{\rho_l-m}.
\end{eqnarray}
Form the first two equations of \eqref{equation}, we get
$$c_{t}=\frac{n_{t}}{m}-\frac{n}{m^{2}}m_{t}=-\frac{mnu_{x}}{m}+\frac{nm^{2}}{m^{2}}u_{x}=0,$$
and
$$\arraycolsep=1.5pt
\begin{array}{rcl}
Q(m)_{t}=\D\left(\frac{m}{\rho_{l}-m}\right)_{t}&=&
\D \left(\frac{1}{\rho_{l}-m}+\frac{m}{(\rho_{l}-m)^{2}}\right)m_{t}\\[3mm]
&=&\D \frac{\rho_{l}}{(\rho_{l}-m)^{2}}m_{t}=-\D\frac{\rho_{l}m^{2}}{(\rho_{l}-m)^{2}}u_{x}\\[3mm]
&=&-\rho_{l}Q(m)^{2}u_{x}.
\end{array}
$$
Then we can rewrite the initial boundary problems
\eqref{equation}-\eqref{initial} into the following forms:
\begin{eqnarray}\label{E2}
\left
\{\begin{array}{l}
 c_t=0\\[3mm]
Q_t+\rho_lQ^2 u_x=0, \\[3mm]
      u_t+(cQ)^\gamma_x=-h(Q)u|u|+g+(c^{\theta}Q^{1+\theta} u_x)_x, \ \  0<x<1, \ \  t>0,\\
  \end{array}
\right.
\end{eqnarray}
with initial data
\begin{eqnarray}
(c Q,u)(x,0)=(c_0Q_0,u_0)(x), \ \ 0\leq x \leq 1,
\end{eqnarray}
and the boundary conditions
\begin{equation}\label{bc}
(cQ)(0,t)=0,\ \ u(1,t)=0,\ \ t>0.
\end{equation}
Here
$$\displaystyle h(Q)=f\rho_l^{2}\frac{Q^2}{(1+Q)^2},\ \ \ f\geq0,
$$
and $c_0,\ Q_0$ is given from  $[n_{0},m_{0}]$ according the
transformations \eqref{tran}. In particular, the first equation of
\eqref{E2} implies that $c(x, t)= c_0(x) := \frac{n_0}{m_0}(x)$.

Let $(cQ_{\infty})(x)$ be the solution of the following stationary
problem:
\begin{equation}
\left \{\begin{array}{l} (cQ_{\infty})^{\gamma}_x=g,\\[2mm]
(cQ_{\infty})(0)=0.
\end{array}
\right.
\end{equation}
Then
\begin{equation}\label{2.6}
(cQ_{\infty})(x)=(gx)^{\frac{1}{\gamma}}.
\end{equation}

Throughout this paper, our assumptions on the initial data $c_0,\
Q_0,\ u_0$ are as follows:\\

$(A_1)$\ \ \ $C_1 x^{\frac{1}{\gamma}}\leq c_0Q_0(x)\leq
 C_2 x^{\frac{1}{\gamma}},\ C_3x^\alpha\leq c_0(x)\leq C_4x^\alpha$, with  some positive constants $0<C_1\leq C_2,\ 0<C_3
 \leq C_4,\  0\leq\alpha<
\dfrac{1}{\gamma}$;
\\
\vspace{1mm}

 $ (A_2)$\ \ \ $u_0\in H^1([0,1])\cap L^5([0,1]),\
 \D\int_0^1x^{\frac{(2+\alpha)\gamma-\theta-1}{\gamma}}((cQ_0)_x^\theta-(cQ_\infty)_x^\theta)^2dx\leq C;$
\\
\vspace{1mm}

 $(A_3)$\ \ \ $(c^{\theta}Q_0^{1+\theta} u_{0x})_x\in L^2([0,1]).$
\\

Under assumptions $(A_1)$-$(A_3)$, we will study the asymptotic
behavior of $cQ$ and $u$ provided that the global weak solution to
the initial boundary value problem \eqref{E2}-\eqref{bc} exists. In
the Lagrangian coordinates, the definition of the weak solution to
\eqref{E2}-\eqref{bc} can be stated as follows:

\bigbreak

\noindent\textbf{Definition 2.1. (Weak solution)}\ \ A pair of
functions $(cQ(x,t),u(x,t))$ is called a weak solution to the
initial boundary problem (\ref{E2})-(\ref{bc}), if
\begin{eqnarray}
c\in C^1([0,1]\times[0,\infty)),\ \ \  cQ, u\in
L^{\infty}([0,1]\times[0,\infty))\cap C^1([0,\infty);L^2(0,1)),
\end{eqnarray}
\begin{eqnarray}
c^\theta Q^{1+\theta}u_x\in L^{\infty}([0,1]\times[0,\infty))\cap
C^{\frac{1}{2}}([0,\infty);L^2(0,1)).
\end{eqnarray}
Furthermore, the following equations hold:
$$
Q_t+\rho_lQ^2u_x=0,\ \ a.e.,
$$
and
$$
\int_0^\infty\int_0^1\left\{u\phi_t+((cQ)^{\gamma}-c^\theta
Q^{1+\theta}
u_x)\phi_x+g\phi-h(Q)u|u|\phi\right\}dxdt+\int_0^1u_0(x)\phi(x,0)dx=0,
$$
for any test functions  $\phi(x,t) \in
C_0^\infty(\Omega)$ with $\Omega= \left\{(x,t): 0\leq x \leq1,\  t\geq 0\right\}.$\\

\noindent\textbf{Remark 2.2.\ \ (Existence of the global weak
solution)}. To our  knowledge, by using the standard line method
(see \cite{FE,Hoff1}
 for example), it is easy to obtain the global existence of the weak solutions to
\eqref{E2}-\eqref{bc}. The details are omitted.

In what follows, $C$  denotes a generic positive constant depending
on initial data, $\gamma$ and $\theta$,  etc., but independent of
$t$. $\delta,\ \epsilon_0$ denote some positive (generally small)
constants.

The main results in this paper can be stated as follows: \bigbreak
\noindent\textbf{Theorem 2.3.}\ \ Let $c_0,\ Q_0,\ u_0$ satisfy
$(A_1)$-$(A_3)$, $\gamma>1$ and
$\theta\in(0,\frac{\gamma}{2}]\cap(0,\gamma-1]\cap(0,1-\alpha\gamma]$.
There exists a constant $0<\epsilon_0<1$, such that if
\begin{eqnarray}
\|u_0\|^2_{L^2}\leq\epsilon_{0},\ \ \ \
\int^1_0x^{1-3/\gamma}(cQ_0-cQ_{\infty})^2dx\leq\epsilon_{0},
\end{eqnarray}
the weak solution $((cQ)(x,t),u(x,t))$ to the initial boundary value
problem (\ref{E2})-(\ref{bc}) satisfies
\begin{eqnarray}\label{cQlimit}
\lim\limits_{t\rightarrow\infty}(cQ)(x,t)=(cQ_{\infty})(x),
\end{eqnarray}
uniformly in $x\in[0,1]$ and
\begin{eqnarray}\label{limitu1}
\lim\limits_{t\rightarrow\infty}\sup\limits_{x\in[\delta,1]}|u(x,t)|=0,
\end{eqnarray}
for any $0<\delta<1$. Furthermore, if we assume further
$\theta<\gamma-1$, then \eqref{limitu1} can be improved as
\begin{eqnarray}\label{limitu2}
\lim\limits_{t\rightarrow\infty}\sup\limits_{x\in[0,1]}|u(x,t)|=0.
\end{eqnarray}

 \bigbreak \noindent\textbf{Theorem 2.4.}\ \ Under the conditions of
Theorem 2.3 and $\theta<\gamma-1$, we have
\begin{eqnarray}\label{2.12}
\|(cQ^{\theta}-cQ^{\theta}_{\infty})(\cdot,t)\|_{L^\infty}\leq
C(1+t)^{-\frac{2\theta}{4\gamma+\alpha\gamma-2}},
\end{eqnarray}
and
\begin{eqnarray}\label{2.13}
\|u(\cdot,t)\|_{L^{\infty}([0,1])}\leq C(1+t)^{-\frac{1}{2}},
\end{eqnarray}
for all $t\geq0$.

\bigbreak

 \noindent\textbf{Remark 2.3.}\ \ If the
coefficient $f=0$ of the frictional force term in \eqref{E2} and
$\alpha=0$, $ C_{3}=C_{4}$ in the assumption $(A_1)$, then
\eqref{E2} is simplified into Navier-Stokes equations with gravity.
In this case, we need not deal with the additional external force
$-h(Q)u|u|$ in \eqref{h(Q)} and derive
$\theta\in(0,\frac{\gamma}{2}]\cap(0,\gamma-1]$. Compared with the
results in \cite{Duan,Fang-Zhang4,Zi-Zhu}, our results can be viewed
as a generalization of the ones in \cite{Duan,Fang-Zhang4,Zi-Zhu}.

\section{Uniform {\it a priori} estimates}
\setcounter{equation}{0}

In this section, we will derive some uniform-in-time  \textit{a
priori} estimates for the solutions to (\ref{E2})-(\ref{bc}) by
classical energy method. First, we list some elementary equalities
which follow from (\ref{E2}) directly. These equalities
will be used frequently later.\\

\noindent\textbf{Lemma 3.1}\ \ Under the conditions of Theorem 2.3,
it holds that for $0<x<1,\  t>0$,
\begin{eqnarray}\label{3.1}
((cQ)^{\theta})_t=\theta\rho_l(cQ_{\infty})^{\gamma}-\theta\rho_l(cQ)^{\gamma}-\theta\rho_l\int^x_0u_tdy-\theta\rho_l\int_0^x h(Q)u|u|dy,
\end{eqnarray}
\begin{eqnarray}\label{3.2}
(cQ)^{\gamma}-(cQ_{\infty})^{\gamma}=-\frac{1}{\theta\rho_l}((cQ)^\theta)_t-\int^x_0u_tdy-\int_0^x h(Q)u|u|dy,
\end{eqnarray}
and
\begin{eqnarray}\label{3.3}
((cQ)^{\theta}-(cQ_{\infty})^{\theta})_{xt}=-\theta\rho_l((cQ)^{\gamma}-(cQ_{\infty})^{\gamma})_x-\theta\rho_l u_t-\theta\rho_l h(Q)u|u|.
\end{eqnarray}

\noindent{\it Proof.} The above three equalities follow from
(\ref{E2}) and (\ref{bc}) directly.
\bigbreak

 The basic energy estimate similar to \cite{Duan} for the case $\theta=1, \gamma=2$
 and \cite{Fang-Zhang4} is given as follows:

\bigbreak \noindent\textbf{Lemma 3.2. (Basic energy estimate).}
Under the conditions of Theorem 2.3, it holds that
\begin{eqnarray}\label{basic energy}
\arraycolsep=1.5pt
\begin{array}[b]{rl}
&\displaystyle\int_0^1\left(\frac{1}{2}u^2+\frac{1}{\rho_l}\int_{Q_{\infty}}^Q
\frac{c^\gamma(h^\gamma-Q_{\infty}^{\gamma})}{h^2}dh\right)dx
+\int_0^t\int_0^1c^{\theta}Q^{1+\theta}u_x^2dxds\\[3mm]
&\displaystyle+\int_0^t\int_0^1h(Q)u^2|u|dxds\leq C\epsilon_{0}, \ \
\ 0<t<\infty,
\end{array}
\end{eqnarray}
where $C$ is a positive constant, independent of $t$ and
$\epsilon_{0}$.

\bigbreak

\vspace{3mm}

\noindent{\it Proof.}  \ Multiplying $(\ref{E2})_2$ and
$(\ref{E2})_3$ by
$\frac{c^{\gamma}(Q^{\gamma}-Q_{\infty}^{\gamma})}{Q^2}$ and $u$
respectively, and integrating the result equations with respect to
$x$ and $t$ over $[0, 1]$ and $[0, t]$,  by using the boundary
conditions \eqref{bc}, we have
$$\arraycolsep=1.5pt
\begin{array}[b]{rl}
&\displaystyle\int_0^1\left(\frac{1}{2}u^2
+\frac{1}{\rho_l}\int_{Q_{\infty}}^Q\frac{c^{\gamma}(h^{\gamma}-Q_{\infty}^{\gamma})}{h^2}dh\right)dx
+\int_0^t\int_0^1h(Q)u^2|u|dxds+\int_0^t\int_0^1c^{\theta}Q^{1+\theta}u_x^2dxds\\[3mm]
=&\displaystyle\int_0^1\left(\frac{1}{2}u_0^2+
\frac{1}{\rho_l}\int_{Q_{\infty}}^{Q_0}\frac{c^{\gamma}(h^{\gamma}-Q_{\infty}^{\gamma})}{h^2}dh\right)dx.
\end{array}
$$
It follows from $(A_1)$ and $(\ref{2.6})$ that \\
$$\arraycolsep=1.5pt
\begin{array}[b]{rcl}
\displaystyle\int_0^1\int_{Q_{\infty}}^{Q_0}\frac{c^{\gamma}(h^{\gamma}-Q_{\infty}^{\gamma})}{h^2}dhdx
&\leq &\displaystyle
C\int_0^1Q_{\infty}^{-2}|Q_0-Q_{\infty}||(cQ_0)^{\gamma}-(cQ_{\infty})^{\gamma}|dx\\
&\leq&\displaystyle C\int_0^1(cQ_{\infty})^{-2}\gamma\xi^{\gamma-1}|cQ_0-cQ_{\infty}|^2dx\\
&\leq &\displaystyle C\int_0^1x^{1-3/\gamma}(cQ_0-cQ_{\infty})^2dx,
\end{array}
$$
where $\xi$ is between $cQ_0$ and $cQ_{\infty}$,\ hence we obtain \\
\begin{eqnarray}
\arraycolsep=1.5pt
\begin{array}[b]{rl}
&\displaystyle\int_0^1\left(\frac{1}{2}u^2+\frac{1}{\rho_l}\int_{Q_{\infty}}^Q
\frac{c^\gamma(h^\gamma-Q_{\infty}^{\gamma})}{h^2}dh\right)dx
+\int_0^t\int_0^1c^{\theta}Q^{1+\theta}u_x^2dxds\\[3mm]
&\displaystyle+\int_0^t\int_0^1h(Q)u^2|u|dxds\leq C\epsilon_{0}.
\end{array}
\end{eqnarray}
This completes the proof.

\vspace{4mm} Combining the the smallness of $\sup_{t\geq0}
\|u(\cdot,t)\|_{L^2(0,1)}$ obtained by the basic energy estimate and
\textit{a priori} assumption, we could get the uniform upper and
lower bounds for $Q(x,t)$ which will play a crucial role in studying
the asymptotic behavior of $cQ$ and $u$.

\bigbreak\noindent\textbf{Lemma 3.3.} Under the conditions in
Theorem 2.3, we have
\begin{eqnarray}\label{cQ}
{K_1}Q_{\infty}(x)\leq Q(x,t)\leq {K_2}Q_{\infty}(x),
\end{eqnarray}
where $(x,t)\in\Omega= \left\{(x,t): 0\leq x \leq1,\  t\geq
0\right\}$, $K_1$ and $K_2$ are two positive constants, independent
of $t$.

\bigbreak\noindent{\it Proof.} Let
$Y(x,t)=(cQ)^{\theta}(x,t)(cQ_{\infty})^{-\theta}(x)$, then due to
(\ref{3.1}), we have
\begin{equation}\label{Yt1}
Y_t=\theta
\rho_l(cQ_{\infty})^{\gamma-\theta}(1-Y^{\frac{\gamma}{\theta}})
-\theta \rho_l\left((cQ_{\infty})^{-\theta}\int_0^x udy\right)_t
-\theta \rho_l(cQ_{\infty})^{-\theta}\int_0^x h(Q)u|u|dy.
\end{equation}
Assume \begin{equation}\label{assumption}
 Y\leq C_5+1,
\end{equation}
 where
$C_5=g^{-\frac{\theta}{\gamma}}C_2^\theta$. The last term on the
right side  of \eqref{Yt1} can be estimated as follows
\begin{equation}\label{h(Q)}
\arraycolsep=1.5pt
\begin{array}[b]{rl}
\D-\theta \rho_l(cQ_{\infty})^{-\theta}\int_0^x h(Q)u|u|dy
&\D\leq C\delta(cQ_\infty)^{-3\theta}\int_0^xh(Q)dy+C(\delta)\int_0^xh(Q)u^2|u|dy\\[3mm]
&\D\leq C\delta(cQ_\infty)^{-3\theta}\int_0^xQ^2dy+C(\delta)\int_0^xh(Q)u^2|u|dy\\[3mm]
&\D\leq C\delta(cQ_\infty)^{-3\theta}\int_0^x(Q^2Q_\infty^{-2})Q_\infty^{2}dy+C(\delta)\int_0^xh(Q)u^2|u|dy\\[3mm]
&\D\leq C\delta(cQ_\infty)^{-3\theta}[\max
Y^{\frac{2}{\theta}}]x^{\frac{2}{\gamma}-2\alpha+1}+C(\delta)\int_0^xh(Q)u^2|u|dy\\[3mm]
 &\D\leq C\delta (C_5+1)^{\frac{2}{\theta}}(cQ_\infty)^{\gamma-\theta}(cQ_\infty)^{-\gamma-2\theta}x^{\frac{2}{\gamma}-2\alpha+1}
 +C(\delta)\int_0^1h(Q)u^2|u|dy\\[3mm]
 &\D\leq C\delta (C_5+1)^{\frac{2}{\theta}}(cQ_\infty)^{\gamma-\theta}+C(\delta)\int_0^1h(Q)u^2|u|dy,
\end{array}
\end{equation}
where we have used Young inequality, and the fact
\begin{eqnarray*}
(cQ_\infty)^{-\gamma-2\theta}x^{\frac{2}{\gamma}-2\alpha+1}\leq C,
\end{eqnarray*}
which follows from \eqref{2.6} and $\D\theta\leq1-\alpha\gamma$.
Then substituting \eqref{h(Q)} into \eqref{Yt1}, we have
$$
\arraycolsep=1.5pt
\begin{array}[b]{rl}
Y_t\leq & \D-\theta \rho_l\left((cQ_{\infty})^{-\theta}\int_0^x
udy\right)_t
+C(\delta)\int_0^1h(Q)u^2|u|dy \\[3mm]
&\displaystyle +\theta
\rho_l(cQ_{\infty})^{\gamma-\theta}(1-Y^{\frac{\gamma}{\theta}}
+C\delta (C_5+1)^{\frac{2}{\theta}}).
\end{array}
$$
Integrating the inequality above over $[0,t]$ with respect to $t$,
we have
\begin{equation}\label{Y}
\arraycolsep=1.5pt
\begin{array}[b]{rl}
Y \leq & \D Y_0+\theta \rho_l(cQ_{\infty})^{-\theta}\int_0^x
(u_0-u)dy
+C(\delta)\int_0^t\int_0^1h(Q)u^2|u|dy \\[3mm]
&\displaystyle +\theta
\rho_l(cQ_{\infty})^{\gamma-\theta}\int_0^t(1-Y^{\frac{\gamma}{\theta}}
+C\delta (C_5+1)^{\frac{2}{\theta}})ds\\[3mm]
 \leq &\D Y_0+2\theta \rho_l(cQ_{_{\infty}})^{-\theta}x^{\frac{1}{2}}\sup\limits_{t\geq0}
 \left(\int^1_0u^2dx\right)^{\frac{1}{2}}
 +C(\delta)\epsilon_{0}+\theta
\rho_l(cQ_{\infty})^{\gamma-\theta}\int_0^t(1-Y^{\frac{\gamma}{\theta}}
+C\delta (C_5+1)^{\frac{2}{\theta}})ds\\[3mm]
\leq &\D
Y_0+C\sup\limits_{t\geq0}\left(\int^1_0u^2dx\right)^{\frac{1}{2}}+C\epsilon_{0}+\theta
\rho_l(cQ_{\infty})^{\gamma-\theta}\int_0^t(1-Y^{\frac{\gamma}{\theta}}
+C \delta(C_{5}+1)^{\frac{2}{\theta}})ds\\[3mm]
\leq &\D Y_0+(C\epsilon_{0})^{\frac{1}{2}}+\theta
\rho_l(cQ_{\infty})^{\gamma-\theta}\int_0^t(1-Y^{\frac{\gamma}{\theta}}
+C\delta (C_5+1)^{\frac{2}{\theta}})ds,
\end{array}
\end{equation}
where we have used the fact
\begin{eqnarray*}
(cQ_\infty)^{-\theta}_{\infty}x^{\frac{1}{2}}\leq C,
\end{eqnarray*}
 which follows from (\ref{2.6}) and the restriction
$\gamma\geq2\theta$ directly. From \eqref{Y}, it is easy to get
\begin{equation}\label{bound of Y}
Y(x,t)\leq\max\{(1+C\delta
(C_5+1)^{\frac{2}{\theta}})^{\frac{\theta}{\gamma}},Y_0+(C\epsilon_{0})^{\frac{1}{2}}\}:=K_2,
\end{equation}
where $Y_0=(cQ_0)^\theta(cQ_\infty)^{-\theta}\leq C_2^\theta
g^{-\frac{\theta}{\gamma}}=C_5$.

 Now, from \eqref{assumption}, \eqref{bound of Y} and continuity of $Y(x,t)$, we can show that
 the \textit{priori} assumption can be closed
provided $\epsilon_{0}$ and $\delta$ sufficiently small. Therefore
the upper bound of $Y(x,t)$ is obtained.

Next, we estimate the lower bound of $Y(x,t)$. Similar to the
estimate of the upper bound of $Y(x,t)$, one has
$$
Y(x,t)\geq Y_0-(C\epsilon_{0})^{\frac{1}{2}}+\theta
\rho_l(cQ_{\infty})^{\gamma-\theta}\int_0^t(1-Y^{\frac{\gamma}{\theta}}
- C \delta)ds.
$$
It is easy to get
$$
Y \ge \min
\{{(1-C\delta)^{\frac{\theta}{\gamma}},Y_0-(C\epsilon_{0})^{\frac{1}{2}}}\}:=K_1.
$$
The proof of Lemma 3.4 is complete.

\bigbreak The next lemma will be
used  for many times later, and it has been proved by Duan in
\cite{Duan} for the case $\theta=1, \gamma=2$.

\bigbreak\noindent\textbf{Lemma 3.4.} Under the conditions in
Theorem 2.3, it holds that
\begin{eqnarray}\label{u2}
\int^t_0\int^1_0u^2dxds\leq C,
\end{eqnarray}
where $C$ is a positive constant, independent of $t$.

\bigbreak\noindent{\it Proof.} Since $u(1,t)=0$, we have
$$
\begin{array}[b]{rl}
|u(x,t)|&\D \leq\int^1_x|u_y(y,t)|dy\\[5mm]
&\D=\int^1_x|u_y(y,t)|c^{\frac {\theta}{2}}Q^{\frac{1+\theta}{2}}c^{-\frac{\theta}{2}}Q^{-\frac{1+\theta}{2}}dy\\[5mm]
&\D\leq\left(\int^1_x
c^{\theta}Q^{1+\theta}u_y^2dy\right)^{\frac{1}{2}}
\left(\int^1_xc^{-\theta}Q^{-(1+\theta)}dy\right)^{\frac{1}{2}}\\[5mm]
&\D\leq C\left(\int^1_x
c^{\theta}Q^{1+\theta}u_y^2dy\right)^{\frac{1}{2}}\left(\int^1_xc^{-(1+\theta)}Q^{-(1+\theta)}dy\right)^{\frac{1}{2}}.
\end{array}
$$
It is easy to see from (\ref{2.6}), (\ref{cQ}) and $(A_{1})$
$$
\begin{array}[b]{rl}
\D\int^1_0\int^1_x(cQ)^{-(1+\theta)}dydx \leq\D
C\int^1_0\int^1_xy^{-\frac{1+\theta}{\gamma}}dydx
\leq\D C, \\[5mm]
\end{array}
$$
provided $\theta<2\gamma-1$, and
\begin{eqnarray}\label{u2esti}
\begin{array}[b]{rl}
&\D\int^1_0u^2dx
\leq\D\int^1_0c^{\theta}Q^{1+\theta}u_x^2dx\int^1_0\int^1_x(cQ)^{-(1+\theta)}dydx
\leq\D C\int^1_0c^{\theta}Q^{1+\theta}u_x^2dx.\\[5mm]
\end{array}
\end{eqnarray}
Using (\ref{basic energy}), we get (\ref{u2}) immediately. This
completes the proof of Lemma 3.4.\\

Next lemma plays a crucial role in dealing with the frictional force
$-h(Q)u|u|$. The proof is based on recurrence method.
 \bigbreak\noindent\textbf{Lemma
3.5.} Under the conditions in Theorem 2.3, it holds that
\begin{eqnarray}\label{pp5}
\int_0^1|u|^5dx+\int_0^t\int_0^1h(Q)|u|^6dxds+\int_0^t\int_0^1c^{\theta}Q^{1+\theta}|u|^3u_x^2dxds
\leq C,
\end{eqnarray}
where $C$ is a positive constant, independent of $t$. Moveover,
\begin{equation}\label{u3u4}
\int_0^t\int_0^1|u|^3dxds\leq C, \ \ \ \
\int_0^t\int_0^1|u|^4dxds\leq C.
\end{equation}
\bigbreak\noindent{\it Proof.} Multiplying $(\ref{E2})_3$ by
$p|u|^{p-2}u$ $(p\geq 2)$ and integrating the result equation over
$[0,1]$ with respect to $x$,
\begin{eqnarray}\label{ddtuP}
\arraycolsep=1.5pt
\begin{array}[b]{rl}
&\displaystyle
\frac{d}{dt}\int_0^1|u|^pdx+p\int_0^1h(Q)|u|^{p-1}u^2dx+p(p-1)\int_0^1c^{\theta}Q^{1+\theta}|u|^{p-2}u_x^2dx\\[5mm]
=&\D p(p-1)\int_0^1|u|^{p-2}u_x(cQ)^{\gamma}dx+gp\int_0^1|u|^{p-2}udx\\[5mm]
\leq& \D\frac{p(p-1)}{2}\int_0^1c^{\theta}Q^{1+\theta}|u|^{p-2}u_x^2dx
+C\int_0^1|u|^{p-2}c^{2\gamma-\theta}Q^{2\gamma-\theta-1}dx
+gp\int_0^1|u|^{p-1}dx\\[5mm]
\leq
&\D\frac{p(p-1)}{2}\int_0^1c^{\theta}Q^{1+\theta}|u|^{p-2}u_x^2dx
+C\int_0^1|u|^{p-2}dx+gp\int_0^1|u|^{p-1}dx,
\end{array}
\end{eqnarray}
 where we have used the fact
\begin{eqnarray*}
c^{2\gamma-\theta}Q^{2\gamma-\theta-1}=c(c Q)^{2\gamma-\theta-1}\leq
C x^{2-\frac{\theta+1}{\gamma}}\leq C.
\end{eqnarray*}
Then integrating \eqref{ddtuP} over $[0,t]$ with respect to $t$, one
has
\begin{equation}\label{p}
\arraycolsep=1.5pt
\begin{array}[b]{rl}
&\displaystyle
\int_0^1|u|^pdx+p\int_0^t\int_0^1h(Q)|u|^{p-1}u^2dxds
+\frac{p(p-1)}{2}\int_0^t\int_0^1c^{\theta}Q^{1+\theta}|u|^{p-2}u_x^2dxds\\[5mm]
&\displaystyle\leq
\int_0^1|u_0|^pdx+C\int_0^t\int_0^1|u|^{p-2}dxds+C\int_0^t\int_0^1|u|^{p-1}dxds.
\end{array}
\end{equation}
Now taking $p=4$ in \eqref{p}, we have
\begin{equation}\label{pp4}
\arraycolsep=1.5pt
\begin{array}[b]{rl}
&\displaystyle
\int_0^1|u|^4dx+4\int_0^t\int_0^1h(Q)|u|^3u^2dxds
+6\int_0^t\int_0^1c^{\theta}Q^{1+\theta}|u|^2u_x^2dxds\\[5mm]
\leq &\displaystyle C\int_0^t\int_0^1|u|^3dxds+C,
\end{array}
\end{equation}
where we have used Lemma 3.4 and the assumption $(A_2)$.

Next we must estimate the first term of the right-hand side in
\eqref{pp4} by using the same method as Lemma 3.4.
 Since $u(1,t)=0$, we have
$$\D
|u|^{\frac{3}{2}}\leq \left|\int_x^1(|u|^{\frac{3}{2}})_ydy\right|
\leq C\int_x^1|u|^{\frac{1}{2}}|u_y|dy \leq
C\left(\int_0^1c^{\theta}Q^{1+\theta}u_y^2dy\right)^{\frac{1}{2}}
\left(\int_x^1c^{-\theta}Q^{-1-\theta}|u|dy\right)^{\frac{1}{2}}.
$$
It is easy to see
\begin{eqnarray}\label{u3}
\arraycolsep=1.5pt
\begin{array}[b]{rl}
\D\int_0^1|u|^3dx
\leq &\displaystyle C\left(\int_0^1c^{\theta}Q^{1+\theta}u_y^2dy\right)\left(\int_0^1\int_x^1c^{-\theta}Q^{-1-\theta}|u|dydx\right)\\[5mm]
\leq &\displaystyle C\left(\int_0^1c^{\theta}Q^{1+\theta}u_y^2dy\right)\left(\int_0^1\int_x^1c(cQ)^{-1-\theta}|u|dydx\right)\\[5mm]
\leq &\displaystyle C\left(\int_0^1c^{\theta}Q^{1+\theta}u_y^2dy\right)\left(\int_0^1\int_x^1y^{-\frac{1+\theta}{\gamma}+\alpha}|u|dydx\right)\\[5mm]
\leq &\displaystyle C\left(\int_0^1c^{\theta}Q^{1+\theta}u_y^2dy\right)\left(\int_0^1\int_x^1
y^{-\frac{4(1+\theta)}{3\gamma}+\frac{4\alpha}{3}}dydx\right)^{\frac{3}{4}}\left(\int_0^1u^4dx\right)^{\frac{1}{4}}\\[5mm]
\leq &\displaystyle
C\int_0^1c^{\theta}Q^{1+\theta}u_y^2dy\left(\int_0^1u^4dx\right)^{\frac{1}{4}},
\end{array}
\end{eqnarray}
where we have used $\theta<(\frac{3}{2}+\alpha)\gamma-1$.

Integrating \eqref{u3} over $[0,t]$ and using \eqref{basic energy},
we get
\begin{eqnarray}\label{intintu3}
\arraycolsep=1.5pt
\begin{array}[b]{rl}
\D\int_{0}^{t}\int_0^1|u|^3dxds
\leq &\displaystyle C\sup_{t\geq0}\left(\int_0^1u^4dx\right)^{\frac{1}{4}}\int_{0}^{t}\int_0^1c^{\theta}Q^{1+\theta}u_y^2dyds\\[5mm]
\leq &\displaystyle
C\epsilon_{0}\sup_{t\geq0}\left(\int_0^1u^4dx\right)^{\frac{1}{4}}.
\end{array}
\end{eqnarray}
Substituting \eqref{intintu3} into \eqref{pp4}, we have

\begin{equation}\label{ppp4}
\arraycolsep=1.5pt
\begin{array}[b]{rl}
&\D\int_0^1|u|^4dx+\int_0^t\int_0^1h(Q)|u|^3u^2dxds+\int_0^t\int_0^1c^\theta
Q^{1+\theta}|u|^2u_x^2dxds\\[3mm]
 \leq &\displaystyle C+C\epsilon_{0}\sup_{t\geq0}\left(\int_0^1u^4dx\right)^{\frac{1}{4}}.
\end{array}
\end{equation}
Young inequality and \eqref{ppp4} imply that
\begin{equation}\label{intu4}
\int_0^1|u|^4dx \leq C.
\end{equation}
Substituting \eqref{intu4} into \eqref{intintu3}, we have
\begin{eqnarray}\label{boundu3}
\int_{0}^{t}\int_0^1|u|^3dxds \leq C.
\end{eqnarray}
Combining \eqref{pp4} and \eqref{boundu3}, we have
\begin{equation}\label{p4}
\D\int_0^1|u|^4dx+\int_0^t\int_0^1h(Q)|u|^3u^2dxds+\int_0^t\int_0^1c^\theta
Q^{1+\theta}|u|^2u_x^2dxds \leq C.
\end{equation}
Taking  $p=5$ in \eqref{p}, we have
\begin{equation}\label{p5}
\arraycolsep=1.5pt
\begin{array}[b]{rl}
&\displaystyle
\int_0^1|u|^5dx+5\int_0^t\int_0^1h(Q)|u|^4u^2dxds
+10\int_0^t\int_0^1c^{\theta}Q^{1+\theta}|u|^3u_x^2dxds\\[5mm]
\leq&\displaystyle
\int_0^1|u_0|^5dx+C\int_0^t\int_0^1|u|^3dxds+C\int_0^t\int_0^1|u|^4dxds.
\end{array}
\end{equation}
 In order to complete the proof of Lemma 3.5, it suffices to
 estimate
the third term of the right-hand side in \eqref{p5} by using the
same method as Lemma 3.4. In fact
$$
|u|^2\leq \left|\int_x^1(|u|^2)_ydy\right|\leq C\int_x^1|u||u_y|dy
\leq
C\left(\int_0^1c^{\theta}Q^{1+\theta}u_y^2dy\right)^{\frac{1}{2}}
\left(\int_x^1c^{-\theta}Q^{-1-\theta}u^2dy\right)^{\frac{1}{2}},
$$
which implies
\begin{eqnarray}\label{u4}
\arraycolsep=1.5pt
\begin{array}[b]{rl}
\D\int_0^1|u|^4dx
\leq &\displaystyle C\left(\int_0^1c^{\theta}Q^{1+\theta}u_y^2dy\right)\left(\int_0^1\int_x^1c^{-\theta}Q^{-1-\theta}u^2dydx\right)\\[5mm]
\leq &\displaystyle C\left(\int_0^1c^{\theta}Q^{1+\theta}u_y^2dy\right)\left(\int_0^1\int_x^1y^{-\frac{(1+\theta)}{\gamma}+\alpha}u^2dydx\right)\\[5mm]
\leq &\displaystyle C\left(\int_0^1c^{\theta}Q^{1+\theta}u_y^2dy\right)\left(\int_0^1\int_x^1
y^{-\frac{5(1+\theta)}{3\gamma}+\frac{5\alpha}{3}}dydx\right)^{\frac{3}{5}}\left(\int_0^1|u|^5dx\right)^{\frac{2}{5}}\\[5mm]
\leq &\displaystyle
C\left(\int_0^1c^{\theta}Q^{1+\theta}u_y^2dy\right)\left(\int_0^1|u|^5dx\right)^{\frac{2}{5}}.
\end{array}
\end{eqnarray}
Here we have used $\theta<(\frac{6}{5}+\alpha)\gamma-1$.

Integrating \eqref{u4} over $[0,t]$ and using \eqref{basic energy},
we get
\begin{eqnarray}\label{intintu4}
\arraycolsep=1.5pt
\begin{array}[b]{rl}
\D\int_{0}^{t}\int_0^1|u|^4dxds
\leq &\displaystyle C\sup_{t\geq0}\left(\int_0^1|u|^5dx\right)^{\frac{2}{5}}\int_{0}^{t}\int_0^1c^{\theta}Q^{1+\theta}u_y^2dyds\\[5mm]
\leq &\displaystyle
C\epsilon_{0}\sup_{t\geq0}\left(\int_0^1|u|^5dx\right)^{\frac{2}{5}}.
\end{array}
\end{eqnarray}
Substituting \eqref{intintu4} into \eqref{p5} and using
\eqref{boundu3} and Young inequality, we have
\begin{equation}\label{intu5}
\int_0^1|u|^5dx \leq C.
\end{equation}
Substituting \eqref{intu5} into \eqref{intintu4}, we have
\begin{eqnarray}\label{boundu4}
\int_{0}^{t}\int_0^1|u|^4dxds \leq C.
\end{eqnarray}
\eqref{p5}, \eqref{boundu3} and \eqref{boundu4} imply \eqref{pp5}.
This completes the proof of Lemma 3.5.

\bigbreak\noindent\textbf{Lemma 3.6.} Under the conditions in
Theorem 2.3, it holds that
\begin{eqnarray}\label{ut2}
\int^1_0u_t^2dx+\int^t_0\int^1_0h(Q)|u|u_t^2dxds+\int^t_0\int^1_0c^{\theta}Q^{1+\theta}u_{xt}^2dxds\leq C,
\end{eqnarray}
where $C$ is a positive constant, independent of $t$.

\bigbreak\noindent{\it Proof.} Differentiating $(\ref{E2})_3$ with
respect to $t$, multiplying $u_t$ and integrating over $[0,1]$ about
$x$, we have
\begin{eqnarray}\label{3.16}
\frac{1}{2}\frac{d}{dt}\int^1_0u_t^2dx+\int^1_0((cQ)^{\gamma})_{xt}u_tdx
=\int^1_0(c^{\theta}Q^{1+\theta}u_x)_{xt}u_tdx-\int^1_0(h(Q)u|u|)_tu_tdx.
\end{eqnarray}
Using the boundary conditions (\ref{bc}), and integrating by parts,
we have
\begin{equation}\label{3.17}
\arraycolsep=1.5pt
\begin{array}[b]{rl}
\displaystyle\int^1_0((cQ)^{\gamma})_{xt}u_tdx
=&\D-\int^1_0((cQ)^{\gamma})_tu_{xt}dx
=-\gamma\int^1_0c^{\gamma}Q^{\gamma-1}Q_tu_{xt}dx
=\D\gamma\rho_l\int^1_0c^{\gamma}Q^{\gamma+1}u_xu_{xt}dx, \\[5mm]
\end{array}
\end{equation}
and
\begin{eqnarray}\label{3.18}
\arraycolsep=1.5pt
\begin{array}[b]{rl}
\displaystyle\int^1_0(c^{\theta}Q^{1+\theta}u_x)_{xt}u_tdx
=&\D-\int^1_0(c^{\theta}Q^{1+\theta}u_x)_tu_{xt}dx\\[5mm]
=&\D-\int^1_0c^{\theta}(Q^{1+\theta})_tu_xu_{xt}dx-\int^1_0c^{\theta}Q^{1+\theta}u_{xt}^2dx\\[5mm]
=&\D(1+\theta)\rho_l\int^1_0c^{\theta}Q^{2+\theta}u_x^2u_{xt}dx-\int^1_0c^{\theta}Q^{1+\theta}u_{xt}^2dx.\\[5mm]
\end{array}
\end{eqnarray}
Substituting (\ref{3.17}) and (\ref{3.18}) into (\ref{3.16}), and using Cauchy
inequality, we get
$$
\arraycolsep=1.5pt
\begin{array}[b]{rl}
&\displaystyle \frac{1}{2}\frac{d}{dt}\int^1_0u_t^2dx+\int^1_0c^{\theta}Q^{1+\theta}u_{xt}^{2}dx
+2\int^1_0h(Q)|u|u_t^2dx\\[5mm]
=&\D
(1+\theta)\rho_l\int^1_0c^{\theta}Q^{2+\theta}u_x^2u_{xt}dx-\gamma\rho_l\int^1_0c^{\gamma}Q^{\gamma+1}u_xu_{xt}dx
+2f\rho_{l}^3\int^1_0\frac{Q^3}{(1+Q)^3}u|u|u_xu_tdx\\[5mm]
\leq&\D\frac{1}{2}\int^1_0c^{\theta}Q^{1+\theta}u_{xt}^2dx+
C\int^1_0c^{\theta}Q^{\theta+3}u_x^4dx
+C\int^1_0c^{2\gamma-\theta}Q^{2\gamma-\theta+1}u_x^2dx
\\[5mm]
&\displaystyle +\int^1_0h(Q)|u|u_t^2dx
+C\int^1_0\frac{Q^6}{h(Q)(1+Q)^6}|u|^3u_x^2dx,
\end{array}
$$
which implies
\begin{eqnarray}\label{3.19}
\arraycolsep=1.5pt
\begin{array}[b]{rl}
&\displaystyle \frac{1}{2}\frac{d}{dt}\int^1_0u_t^2dx+\frac{1}{2}\int^1_0c^{\theta}Q^{1+\theta}u_{xt}^{2}dx
+\int^1_0h(Q)|u|u_t^2dx\\[5mm]
\leq & \D C\int^1_0c^{2\gamma-\theta}Q^{2\gamma-\theta+1}u_x^2dx
+C\int^1_0c^{\theta}Q^{\theta+3}u_x^4dx
+C\int^1_0\frac{Q^4}{(1+Q)^4}|u|^3u_x^2dx\\[5mm]
\leq &\displaystyle
C\int^1_0(c^{2\gamma-2\theta}Q^{2\gamma-2\theta})c^{\theta}Q^{1+\theta}u_x^2dx
+C\int^1_0c^{\theta}Q^{1+\theta}u_x^2(Qu_x)^2dx
+C\int^1_0 c^{-\theta}Q^{3-\theta}c^{\theta}Q^{1+\theta}|u|^3u_x^2dx\\[5mm]
\leq &\displaystyle C \int^1_0c^{\theta}Q^{1+\theta}u_x^2dx +C
\int^1_0c^{\theta}Q^{1+\theta}u_x^2(Qu_x)^2dx
+C\int^1_0c^{\theta}Q^{1+\theta}|u|^3u_x^2dx,
\end{array}
\end{eqnarray}
where we have used the fact
$$
c^{-\theta}Q^{3-\theta}\leq C x^{\frac{3(1-\alpha
\gamma)}{\gamma}-\frac{\theta}{\gamma}}\leq C,
$$
which follows from $\theta\leq 1-\alpha \gamma$ and \eqref{2.6}. By
\eqref{2.6} and Lemma 3.5, we get
\begin{equation}\label{3.20}
\arraycolsep=1.5pt
\begin{array}[b]{rcl}
(Qu_x)^2&=&\D(c^{\theta}Q^{1+\theta}u_x)^2(cQ)^{-2\theta}\\[5mm]
&=&\D(cQ)^{-2\theta}\left(\int^x_0u_tdy+(cQ)^{\gamma}-(cQ_{\infty})^{\gamma}+\int_0^xh(Q)u|u|dy\right)^2\\[5mm]
&\leq&\D
C(cQ)^{-2\theta}\left\{\left(\int^x_0u_tdy\right)^2+(cQ)^{2\gamma}+(cQ_{\infty})^{2\gamma}
+\left(\int_0^xh(Q)u|u|dy\right)^2\right\}\\[5mm]
&\leq&\D C(cQ)^{-2\theta}x\int^1_0u_t^2dx+C(cQ)^{2\gamma-2\theta}+C(cQ)^{-2\theta}x^2+C(cQ)^{-2\theta}x\int^1_0u^4dx\\[5mm]
&\leq&\D C\int^1_0u_t^2dx+C.\\[5mm]
\end{array}
\end{equation}
Here we have used the fact
\begin{equation*}
(cQ)^{-2\theta}x\leq C x^{1-\frac{2\theta}{\gamma}} \leq C,
\end{equation*}
which follows from $\theta\leq\dfrac{\gamma}{2}$.
 Then it follows from (\ref{3.19}), (\ref{3.20}) and Lemma 3.5 that
$$
\arraycolsep=1.5pt
\begin{array}[b]{rl}
 &\displaystyle
 \frac{d}{dt}\int^1_0u_t^2dx+\int^1_0c^{\theta}Q^{1+\theta}u_{xt}^2dx+2\int^1_0h(Q)|u|u_t^2dx\\[3mm]
\leq & C \displaystyle
 \int^1_0c^{\theta}Q^{1+\theta}u_x^2dx\int^1_0u^2_tdx+C\int^1_0c^{\theta}Q^{1+\theta}u_x^2dx.
\end{array}
$$
According to assumption $(A_3)$, and taking care of $(\ref{E2})_3$,
we obtain $\displaystyle\int^1_0u_{t}^2(x,0)dx\leq C$. Then Gronwall
inequality and (\ref{basic energy}) imply that
$$
\int^1_0u_t^2dx\leq
C+C\exp\left(C\int^t_0\int^1_0c^{\theta}Q^{1+\theta}u_x^2dxds\right)\leq C,
$$
and (\ref{ut2}) follows. This proves Lemma 3.6.

\bigbreak\noindent\textbf{Corollary 3.7.} It follows from (\ref{ut2}) and
(\ref{3.20}) that
\begin{eqnarray}
\left\|Q u_x\right\|_{L^{\infty}([0,1]\times[0,\infty))}\leq C,
\end{eqnarray}
where $C$ is a positive constant, independent of $t$.

\bigbreak\noindent\textbf{Lemma 3.8.} Assume the conditions in
Theorem
2.3 hold, and let $\beta$ be a fixed constant.\\
(i) if $\beta>0$, then
\begin{eqnarray}\label{3.22}
\int^t_0\int^1_0((cQ)^{\gamma}-(cQ_{\infty})^{\gamma})((cQ)^{\beta}-(cQ_{\infty})^{\beta})dxds\leq
C,
\end{eqnarray}
(ii) if $\beta\in W$, then
\begin{eqnarray}\label{3.23}
\int^t_0\int^1_0((cQ)^{\gamma}-(cQ_{\infty})^{\gamma})((cQ_{\infty})^{-\beta}-(cQ)^{-\beta})dxds\leq
C,
\end{eqnarray}
where $C$ is a positive constant, independent of $t$.
$W=(0,2(1-\gamma\alpha)]\cap(0,\frac{(1-\alpha)\gamma+1-\theta}{2}]$.

\bigbreak\noindent{\it Proof.} We only prove (ii) here, since the
proof of (i) is almost the same with the later.

 Due to (\ref{3.2}), we
have
$$
\arraycolsep=1.5pt
\begin{array}[b]{rl}
&\D\int^t_0\int^1_0((cQ)^{\gamma}-(cQ_{\infty})^{\gamma})((cQ_{\infty})^{-\beta}-(cQ)^{-\beta})dxds\\[5mm]
=&\D-\frac{1}{\theta\rho_l}\int^t_0\int^1_0((cQ)^{\theta})_t((cQ_{\infty})^{-\beta}-(cQ)^{-\beta})dxds\\[5mm]
&\D -\int^t_0\int^1_0
((cQ_{\infty})^{-\beta}-(cQ)^{-\beta})\left(\int^x_0udy\right)_tdxds\\[5mm]
&\displaystyle
-\int^t_0\int^1_0\left(\int^x_0h(Q)u|u|dy\right)((cQ_{\infty})^{-\beta}-(cQ)^{-\beta})dxds
\\[5mm]
= &\D I_1+I_2+I_3,\\[5mm]
\end{array}
$$
where
$$
\arraycolsep=1.5pt
\begin{array}[b]{rl}
I_1&\D=-\frac{1}{\rho_l\theta}\int^t_0\int^1_0((cQ)^{\theta})_t((cQ_{\infty})^{-\beta}-(cQ)^{-\beta})dxds\\[5mm]
&\D=\frac{1}{\rho_l\theta}\int^1_0((cQ_0)^{\theta}-(cQ)^{\theta})(cQ_{\infty})^{-\beta}dx
+\frac{1}{\rho_l}\int^1_0\int^t_0(cQ)^{\theta-\beta-1}(cQ)_tdsdx\\[5mm]
&\D=I_1^{(1)}+I_1^{(2)}.
\end{array}
$$
By using $(A_1)$, (\ref{2.6}) and (\ref{cQ}), one easily gets
\begin{eqnarray}\label{3.24}
\left|I_1^{(1)}\right|\leq
C\int^1_0x^{\frac{\theta-\beta}{\gamma}}dx\leq C,
\end{eqnarray}
provided $\beta<\theta+\gamma$.

Next we estimate $I_1^{(2)}$ by two cases.

\bigbreak \noindent{\bf Case 1:} $ \beta\neq\theta.$

In this case, if $\beta<\theta+\gamma$, similarly to (\ref{3.24}),
we have
$$
\left|I_1^{(2)}\right|=\left|\frac{1}{\rho_l(\theta-\beta)}\int^1_0((cQ)^{\theta-\beta}-(cQ_{0})^{\theta-\beta})dx\right|
\leq C.
$$
\bigbreak \noindent{\bf Case 2:} $\beta=\theta.$

In this case, also by $(A)_1$, (\ref{cQ}), we have
$$
\left|I_1^{(2)}\right|=\left|\frac{1}{\rho_{l}}\int^1_0(\ln(cQ)-\ln(cQ_0))dx\right|
\leq C\int^1_0|\ln x|dx \leq C.
$$
Hence, if $\beta<\theta+\gamma$, we have
\begin{eqnarray}\label{3.25}
\left|I_1\right|\leq C.
\end{eqnarray}

Now we estimate $I_2$ as follows:
$$
\arraycolsep=1.5pt
\begin{array}[b]{rl}
I_2=&\D-\int^1_0((cQ_{\infty})^{-\beta}-(cQ)^{-\beta})\left(\int^x_0udy\right)dx
+\int^1_0((cQ_{\infty})^{-\beta}-(cQ_{0})^{-\beta})\left(\int^x_0u_0dy\right)dx\\[5mm]
&\D+\int^t_0\int^1_0((cQ_{\infty})^{-\beta}-(cQ)^{-\beta})_t\left(\int^x_0udy\right)dxds\\[5mm]
=& I_2^{(1)}+I_2^{(2)},\\[5mm]
\end{array}
$$
where
\begin{equation}\label{3.26}
\arraycolsep=1.5pt
\begin{array}[b]{rl}
\left|I_2^{(1)}\right|&\D=\left|-\int^1_0((cQ_{\infty})^{-\beta}-(cQ)^{-\beta})\left(\int^x_0udy\right)dx
+\int^1_0((cQ_{\infty})^{-\beta}-(cQ_{0})^{-\beta})\left(\int^x_0u_0dy\right)dx\right|\\[5mm]
&\D\leq
C\int^1_0x^{-\frac{\beta}{\gamma}}x^{\frac{1}{2}}\left\{\left(\int^1_0u^2dx\right)^{\frac{1}{2}}
+\left(\int^1_0u^2_0dx\right)^{\frac{1}{2}}\right\}dx\\[5mm]
&\D\leq C\int^1_0x^{\frac{1}{2}-\frac{\beta}{\gamma}}dx\leq C.
\end{array}
\end{equation}
Here, $\beta<\frac{3}{2}\gamma$, $(A_1)$, (\ref{2.6}), (\ref{basic
energy}), (\ref{cQ}) and H\"{o}lder inequality were used.

 From
$(\ref{E2})_2$, H\"{o}lder inequality, (\ref{basic energy}),
(\ref{cQ}) and (\ref{u2}), we have
\begin{eqnarray}\label{3.27}
\arraycolsep=1.5pt
\begin{array}[b]{rl}
\left|I_2^{(2)}\right|
&\displaystyle =\int^t_0\int^1_0((cQ_{\infty})^{-\beta}-(cQ)^{-\beta})_t\left(\int^x_0udy\right)dxds\\[5mm]
&\D\leq C\int^t_0\int^1_0c^{-\beta}Q^{1-\beta}\left|u_x\right|x^{\frac{1}{2}}
\left(\int^x_0u^2dy\right)^{\frac{1}{2}}dxds\\[5mm]
&\D\leq C\left(\int^t_0\int^1_0xc^{-\theta-2\beta}Q^{1-\theta-2\beta}c^{\theta}Q^{1+\theta}u_x^2dxds\right)^{\frac{1}{2}}
\left(\int^t_0\int^1_0u^2dxds\right)^{\frac{1}{2}}\\[5mm]
&\D\leq
C\sup\limits_{[0,1]\times[0,\infty)}(x^{1-\alpha}(cQ)^{1-\theta-2\beta})^{\frac{1}{2}}
\left(\int^t_0\int^1_0c^{\theta}Q^{1+\theta}u_x^2dxds\right)^{\frac{1}{2}}\\[5mm]
&\D\leq C,
\end{array}
\end{eqnarray}
provided $\beta\leq\frac{(1-\alpha)\gamma-\theta+1}{2}$.

 Next, we
estimate $I_3$,
\begin{equation}\label{3.28}
\arraycolsep=1.5pt
\begin{array}[b]{rl}
I_3
&\displaystyle=\int^t_0\int^1_0\left(\int^x_0h(Q)u|u|dy\right)((cQ_{\infty})^{-\beta}-(cQ)^{-\beta})dxds\\[5mm]
&\displaystyle \leq
\int^t_0\int^1_0\left(\int^x_0h(Q)u|u|dy\right)\left\{(cQ_{\infty})^{-\beta}+(cQ)^{-\beta}\right\}dxds\\[5mm]
&\displaystyle = I_3^{(1)}+I_3^{(2)}.
\end{array}
\end{equation}
Because of the similarity of $I_3^{(1)}$ and $I_3^{(2)}$, we only
estimate one of them,
\begin{equation}\label{3.29}
\left|I_3^{(1)}\right| \leq
\sup\left[(cQ_{\infty})^{-\beta}x^{\frac{2}{\gamma}-2\alpha}\right]\int^t_0\int^1_0u^2
dxds\leq C,
\end{equation}
provided $\beta\leq2-2\gamma\alpha$.\\[3mm]

$\gamma>1$ and
$\theta\in(0,\frac{\gamma}{2}]\cap(0,\gamma-1]\cap(0,1-\alpha\gamma]$
implies that
$1-\alpha\gamma\in(0,2(1-\gamma\alpha)]\cap(0,\frac{(1-\alpha)\gamma+1-\theta}{2}]
\cap(0,\gamma+\theta)\cap(0,\frac{3}{2}\gamma)=(0,2(1-\gamma\alpha)]\cap(0,\frac{(1-\alpha)\gamma+1-\theta}{2}]$.
Hence, $W$ is not empty, it follows from (\ref{3.24})-(\ref{3.29})
that (\ref{3.23}) holds. This completes the proof of Lemma 3.8.

In the next lemma, we give the uniform estimate of the derivative of
the function $(cQ)(x,t)$ in a appropriate weighted $L^2$ space which
is a little different from \cite{Fang-Zhang4,Zi-Zhu}. Similar
estimate for Navier-Stokes equations has been established by Zhang
and Fang in \cite{Fang-Zhang4} for
$\theta\in(0,\gamma-1)\cap(0,\frac{\gamma}{2}]$, by Zhu and Zi in
\cite{Zi-Zhu} for $\theta\in(0,\gamma-1]\cap(0,\frac{\gamma}{2}]$.

 \bigbreak
\noindent\textbf{Lemma 3.9.} Assume the conditions in Theorem 2.3
hold, then
\begin{equation}\label{3.30}
\int^1_0x^{\frac{(2+\alpha)\gamma-\theta-1}{\gamma}}((cQ)^{\theta}-(cQ_{\infty})^{\theta})^2_xdx
+\int^t_0\int^1_0x^{\frac{(3+\alpha)\gamma-2\theta-1}{\gamma}}((cQ)^{\theta}-(cQ_{\infty})^{\theta})^2_xdxdt\leq
C,
\end{equation}
where $C$ is a positive constant, independent of $t$. \bigbreak
\noindent{\it Proof.} From (\ref{3.3}) and (\ref{2.6}), we get
\begin{eqnarray}\label{3.31}
\arraycolsep=1.5pt
\begin{array}[b]{rl}
&\left(((cQ)^{\theta})_x-((cQ_{\infty})^{\theta})_x+\theta\rho_l
u\right)_t+\gamma\rho_l(cQ)^{\gamma-\theta}\left(((cQ)^{\theta})_x-((cQ_{\infty})^{\theta})_x+\theta\rho_l
u\right)\\[5mm]
=&\theta\rho_l^2\gamma(cQ)^{\gamma-\theta}u+\theta\rho_l
g\left(1-\frac{(cQ)^{\gamma-\theta}}{(cQ_{\infty})^{\gamma-\theta}}\right)-\theta\rho_l
h(Q)u|u|.
\end{array}
\end{eqnarray}
Multiplying (\ref{3.31}) by
$x^{\frac{(2+\alpha)\gamma-\theta-1}{\gamma}}\left(((cQ)^\theta)_x-((cQ_{\infty})^{\theta})_x+\theta\rho_l
u\right)$, and integrating the result equation on $[0, 1]\times[0,
t]$, we get
\begin{equation*}
\arraycolsep=1.5pt
\begin{array}[b]{rl}
&\displaystyle
\frac{1}{2}\int^1_0x^{\frac{(2+\alpha)\gamma-\theta-1}{\gamma}}\left(((cQ)^{\theta})_x-((cQ_{\infty})^{\theta})_x+\theta\rho_l
u\right)^2dx\\[5mm]
&\displaystyle+\gamma\rho_l\int^t_0\int^1_0(cQ)^{\gamma-\theta}x^{\frac{(2+\alpha)\gamma-\theta-1}{\gamma}}
\left(((cQ)^{\theta})_x-((cQ_{\infty})^{\theta})_x+\theta\rho_l
u\right)^2dxds\\[5mm]
=&\displaystyle
\frac{1}{2}\int^1_0x^{\frac{(2+\alpha)\gamma-\theta-1}{\gamma}}\left(((cQ_0)^{\theta})_x-((cQ_{\infty})^{\theta})_x+\theta\rho_l
u_0\right)^2dx\\[5mm]
&\displaystyle
+\int^t_0\int^1_0x^{\frac{(2+\alpha)\gamma-\theta-1}{\gamma}}\left(((cQ)^{\theta})_x-((cQ_{\infty})^{\theta})_x+\theta\rho_l
u\right)\left(\theta\rho_l^2\gamma(cQ)^{\gamma-\theta}u+\theta\rho_l
g\left(1-\frac{(cQ)^{\gamma-\theta}}{(cQ_{\infty})^{\gamma-\theta}}\right)\right)dxds\\[5mm]
&\displaystyle+\int^t_0\int^1_0x^{\frac{(2+\alpha)\gamma-\theta-1}{\gamma}}
\left(((cQ)^\theta)_x-((cQ_{\infty})^{\theta})_x+\theta\rho_{l}
u\right)(-\theta\rho_l h(Q)u|u|)dx d s
\end{array}
\end{equation*}
\begin{equation}\label{3.32}
\arraycolsep=1.5pt
\begin{array}[b]{rl}
 \leq &\displaystyle
C+\gamma\rho_{l}\frac{1}{2}\int^t_0\int^1_0(cQ)^{\gamma-\theta}x^{\frac{(2+\alpha)\gamma-\theta-1}{\gamma}}
\left(((cQ)^{\theta})_x-((cQ_{\infty})^{\theta})_x+\theta\rho_l
u\right)^2dx d s\\[5mm]
&\displaystyle
+C\int^t_0\int^1_0(cQ)^{\gamma-\theta}x^{\frac{(2+\alpha)\gamma-\theta-1}{\gamma}}u^2dxds
+
C\int^t_0\int^1_0x^{\frac{(2+\alpha)\gamma-\theta-1}{\gamma}}(cQ)^{\theta-\gamma}
\left(1-\frac{(cQ)^{\gamma-\theta}}{(cQ_{\infty})^{\gamma-\theta}}\right)^2dx d s\\[5mm]
&\displaystyle+C\int^t_0\int^1_0x^{\frac{(2+\alpha)\gamma-\theta-1}{\gamma}}(cQ)^{\theta-\gamma}h^2(Q)u^4dx d s\\[5mm]
\leq&\displaystyle
C+\frac{1}{2}\int^t_0\int^1_0\gamma\rho_{l}(cQ)^{\gamma-\theta}x^{\frac{(2+\alpha)\gamma-\theta-1}{\gamma}}
\left(((cQ)^{\theta})_x-((cQ_{\infty})^{\theta})_x+\theta\rho_l
u\right)^2dxds\\[5mm]
&\displaystyle+C\int^t_0\int^1_0 x^{\frac{(3+\alpha)\gamma-2\theta-1}{\gamma}}u^2dxds
+C\int^t_0\int^1_0x^{-\frac{(1-\alpha)\gamma+1}{\gamma}}\left((cQ)^{\gamma}-(cQ_{\infty})^{\gamma}\right)^2dxds\\[5mm]
&\displaystyle
+C\int^t_0\int^1_0x^{\frac{(1+\alpha)\gamma-1}{\gamma}}u^4dxds,
\end{array}
\end{equation}
which implies

\begin{equation}\label{3.35}
\arraycolsep=1.5pt
\begin{array}[b]{rl}
&\displaystyle
\int^1_0x^{\frac{(2+\alpha)\gamma-\theta-1}{\gamma}}\left(((cQ)^{\theta})_x-((cQ_{\infty})^{\theta})_x+\theta\rho_l
u\right)^2dx\\[5mm]
&\displaystyle+\gamma\rho_l\int^t_0\int^1_0(cQ)^{\gamma-\theta}x^{\frac{(2+\alpha)\gamma-\theta-1}{\gamma}}
\left(((cQ)^{\theta})_x-((cQ_{\infty})^{\theta})_x+\theta\rho_l
u\right)^2dxds\\[5mm]
\leq & \displaystyle C+C\int^t_0\int^1_0
x^{\frac{(3+\alpha)\gamma-2\theta-1}{\gamma}}u^2dxds+C\int^t_0\int^1_0x^{\frac{(1+\alpha)\gamma-1}{\gamma}}u^4dxds\\[5mm]
&+\displaystyle \int^t_0\int^1_0
A(x,t)\left\{(cQ)^{\gamma}-(cQ_{\infty})^{\gamma}\right\}\left\{(cQ)^{-\beta}-(cQ_{\infty})^{-\beta}\right\}dxds,\\[3mm]
\end{array}
\end{equation}
where
\begin{equation*}
\arraycolsep=1.5pt
\begin{array}[b]{rl}
A(x,t)=&\displaystyle\frac{x^{\frac{(2+\alpha)\gamma-\theta-1}{\gamma}}(cQ)^{\theta-\gamma}(cQ_{\infty})^{2\theta-2\gamma}\left\{
(cQ_{\infty})^{\gamma-\theta}-(cQ)^{\gamma-\theta}\right\}^2}{\left\{(cQ)^{\gamma}-
(cQ_{\infty})^{\gamma}\right\}\left\{(cQ_{\infty})^{-\beta}-
(cQ)^{-\beta}\right\}}\\[5mm]
=&\displaystyle\frac{x^{\frac{(2+\alpha)\gamma-\theta-1}{\gamma}}(cQ)^{\beta+\theta-\gamma}(cQ_{\infty})^{\beta+2\theta-2\gamma}\left\{
(cQ_{\infty})^{\gamma-\theta}-(cQ)^{\gamma-\theta}\right\}^2}{\left\{(cQ)^{\gamma}-
(cQ_{\infty})^{\gamma}\right\}\left\{(cQ)^{\beta}-
(cQ_{\infty})^{\beta}\right\}}.\\[5mm]
\end{array}
\end{equation*}
It is easy to see that
\begin{equation*}
\arraycolsep=1.5pt
\begin{array}[b]{rl}
|A(x,t)|\leq  & C
x^{\frac{(2+\alpha)\gamma-\theta-1}{\gamma}}(cQ_{\infty})^{\beta-2\gamma+\theta}\\[5mm]
\leq & C
x{\frac{(2+\alpha)\gamma-\theta-1+\beta-2\gamma+\theta}{\gamma}}\\[5mm]
\leq & C,
\end{array}
\end{equation*}
provided $\beta=1-\alpha\gamma$. Here we have used Lemma 3.4 and
\eqref{u3u4}. Taking $\beta=1-\alpha\gamma$ in (\ref{3.23}), and
using Lemma 3.8, we obtain (\ref{3.30}) immediately from
(\ref{3.35}).

This completes the proof of Lemma 3.9.

\section{Asymptotic behavior}
\setcounter{equation}{0} In this section, based on the \emph{a
priori} estimates obtained in Section 3, we will consider the
asymptotic behavior of the solution $((cQ)(x,t),u(x,t))$ to the
initial boundary problem (\ref{E2})-(\ref{bc}). We will show both
$(cQ)(x,t)$ and $u(x,t)$ converge to the stationary state uniformly
as $t\rightarrow\infty.$

To apply the uniform estimates obtained above to study the
asymptotic behavior of $(cQ)(x,t)$ and $u(x,t)$, we introduce the
following lemma (cf. \cite{S-Z,Fang-Zhang4, Duan}), and omit the
details of the proof. \bigbreak \noindent\textbf{Lemma 4.1.} Suppose
that $y\in W^{1,1}_{loc}(\mathbb{R^{+}})$ satisfies
$$
y=y_1'+y_2,
$$
and
$$
\left|y_2\right|\leq\sum^n_{i=1}\alpha_i,\ \
\left|y'\right|\leq\sum^n_{i=1}\beta_i,\ \ \textrm{on}\ \
\mathbb{R^{+}},
$$
where $y_1\in
W^{1,1}_{loc}(\mathbb{R^+}),\lim\limits_{s\rightarrow+\infty}y_1(s)=0$
and $\alpha_i,\ \beta_i\in L^{p_i}(\mathbb{R^+})$ for some
$p_i\in[1,\infty),\  i=1,\cdots, n.$ Then
$$
\lim\limits_{s\rightarrow+\infty}y(s)=0.
$$

{\noindent\underline{\it Proof of Theorem 2.3}}.

At first, we consider the convergence of $u(x,t)$. To this end,
motivated by \cite{Okada1}, we introduce a function of $t$ as
follows:
$$
f(t)=\int^1_0c^{\theta}Q^{1+\theta}_{\infty}u_x^2dx.
$$
Then by Lemma 3.2 and Lemma 3.3, we get
$$
\int_0^\infty|
f(t)|dt=\int_0^\infty\int_0^1c^{\theta}Q_{\infty}^{1+\theta}u_x^2dxdt\leq
C.
$$
Furthermore, by Lemmas 3.2, 3.3, 3.6 and Cauchy inequality, we have
$$\int_0^\infty\left|\frac{df(t)}{dt}\right|dt
\leq
\int_0^\infty\int_0^1c^{\theta}Q^{1+\theta}_{\infty}u_x^2dxdt+\int_0^\infty\int_0^1c^{\theta}Q^{1+\theta}_{\infty}u_{xt}^2dxdt
\leq C.$$
 It follows from Lemma 4.1 that
$$\lim_{t\rightarrow\infty}\int_0^1c^{\theta}Q^{1+\theta}_{\infty}u_x^2dx=0.$$
Therefore, by H\"{o}lder inequality, we have for $x\in[\delta,1]$
with $0<\delta<1$,
\begin{eqnarray}\label{4.1}
\arraycolsep=1.5pt  \begin{array}[b]{rl}
 \mid
u(x,t)\mid &\D =\left|\int_x^1u_z(z,t)dz
\right|=\left|\int_x^1c^{\frac{\theta}{2}}Q^{\frac{1+\theta}{2}}_{\infty}u_z(z,t)
c^{-\frac{\theta}{2}}Q^{-\frac{1+\theta}{2}}_{\infty}dz\right|\\[5mm]
&\D \leq
\left(\int_0^1c^{\theta}Q^{1+\theta}_{\infty}u_x^2(x,t)dx\right)^\frac{1}{2}
\left(\int_{\delta}^1c^{-\theta}Q^{-(1+\theta)}_{\infty}dx\right)^\frac{1}{2}\\[5mm]
 &\D \leq C\left(\int_{0}^1c^{\theta}Q^{1+\theta}_{\infty}u_x^2(x,t)dx\right)^\frac{1}{2}\rightarrow
0,\\[5mm]
\end{array}
\end{eqnarray}
as $t\rightarrow \infty$, i.e.,
\begin{eqnarray}
\lim_{t\rightarrow\infty}\sup_{x\in[\delta,1]}u(x,t)=0.
\end{eqnarray}

Furthermore, when $\delta=0$,
\begin{eqnarray*}
\int_{0}^1c^{-\theta}Q^{-(1+\theta)}_{\infty}dx \leq
\int_{0}^1x^{-\frac{1+\theta}{\gamma}+\alpha}dx \leq C,
\end{eqnarray*}
provided $\theta<\gamma-1$. Thus, when $\theta<\gamma-1$, we have
\begin{eqnarray}
\lim_{t\rightarrow\infty}\sup_{x\in[0,1]}u(x,t)=0.
\end{eqnarray}

 Next we consider the convergence of
$(cQ)(x,t)$. Firstly, we shall show that $(cQ)(x,t)$ tends to the
stationary state $(cQ_\infty)(x)$ in the sense of integral as
$t\rightarrow\infty$. The similar conclusion for Navier-Stokes
equation was obtained in \cite{Fang-Zhang4} and \cite{Duan} before.

Taking $\beta=\gamma$ in (\ref{3.22}), we have
\begin{eqnarray}\label{4.3}
\int^1_0((cQ)^{\gamma}-(cQ_{\infty})^{\gamma})^2dx\in
L^1(\mathbb{R^+}).
\end{eqnarray}
On the other hand, it is easy to see that
\begin{eqnarray}\label{4.4}
\arraycolsep=1.5pt
\begin{array}[b]{rl}
\D\left|\frac{d}{dt}\int^1_0((cQ)^{\gamma}-(cQ_{\infty})^{\gamma})^2dx\right|
=&\D\left|\int^1_02((cQ)^{\gamma}-(cQ_{\infty})^{\gamma})\gamma c^{\gamma}Q^{\gamma-1}Q_tdx\right|\\[5mm]
=&\D2\gamma\rho_l\left|\int^1_0((cQ)^{\gamma}-(cQ_{\infty})^{\gamma})c^{\gamma}Q^{\gamma+1}u_xdx\right|\\[5mm]
\leq&\D
C\left(\int^1_0c^{2\gamma-\theta}Q^{2\gamma+1-\theta}dx\right)^{1/2}\left(\int^1_0c^{\theta}Q^{1+\theta}u_x^2dx\right)^{1/2}\\[5mm]
\leq&\D C\left(\int^1_0c^{2\gamma-\theta}Q^{2\gamma-\theta}dx\right)^{1/2}\left(\int^1_0c^{\theta}Q^{1+\theta}u_x^2dx\right)^{1/2}\\[5mm]
\leq&\D C\left(\int^1_0c^{\theta}Q^{1+\theta}u_x^2dx\right)^{1/2}.\\[5mm]
\end{array}
\end{eqnarray}
Then \eqref{4.3}, \eqref{4.4} and Lemma 4.1 yield
\begin{eqnarray}
\lim_{t\rightarrow\infty}\int^1_0((cQ)^{\gamma}-(cQ_{\infty})^{\gamma})^2dx=0.
\end{eqnarray}
Consequently,
\begin{eqnarray}\label{4.6}
\arraycolsep=1.5pt
\begin{array}[b]{rl}
\D\int^1_0(cQ-cQ_{\infty})^{2\gamma}dx
=&\D\int^1_0\frac{(cQ-cQ_{\infty})^{2\gamma}}{((cQ)^{\gamma}-(cQ_{\infty})^{\gamma})^2}
((cQ)^{\gamma}-(cQ_{\infty})^{\gamma})^2dx\\[5mm]
\leq&\D
C\int^1_0((cQ)^{\gamma}-(cQ_{\infty})^{\gamma})^2dx\rightarrow0,\ \ \textrm{as}\ \  t\rightarrow\infty.\\[5mm]
\end{array}
\end{eqnarray}
 For $q\in(0,2\gamma)$,  we have from \eqref{4.6} and H\"{o}lder  inequality  that
\begin{eqnarray}
\begin{array}[b]{rl}\label{4.8}
&\D\int^1_0\left|cQ-cQ_{\infty}\right|^{q}dx
\leq\D C\left(\int^1_0(cQ-cQ_{\infty})^{2\gamma}dx\right)^{\frac{q}{2\gamma}}\rightarrow0.\\[5mm]
\end{array}
\end{eqnarray}
On the other hand, for $q\in(2\gamma,\infty)$, we have from (\ref{2.6})
and (\ref{cQ}) that
\begin{eqnarray}\label{4.9}
\arraycolsep=1.5pt
\begin{array}[b]{rl}
\D\int^1_0\left|cQ-cQ_{\infty}\right|^{q}dx
=&\D\int^1_0\left|cQ-cQ_{\infty}\right|^{q-2\gamma}(cQ-cQ_{\infty})^{2\gamma}dx\\[5mm]
\leq&\D C\int^1_0(cQ-cQ_{\infty})^{2\gamma}dx\rightarrow0.\\[5mm]
\end{array}
\end{eqnarray}
Hence, by \eqref{4.8} and \eqref{4.9}, we get
\begin{eqnarray}\label{CQLq}
\|(cQ-cQ_{\infty})(\cdot,t)\|_{L^q}\rightarrow0, \ \ \textrm{as}\ \
t\rightarrow\infty,\ \ q\in(0,\infty).
\end{eqnarray}
We are now in a position to show the uniform convergence of
$(cQ)(x,t)$. To this end, choosing a positive number $k$ large
enough, which is to be determined later, applying H\"{o}lder
inequality, (\ref{2.6}) and (\ref{cQ}), we have by
$(cQ)(0,t)=(cQ_{\infty})(0)=0$
\begin{eqnarray}\label{4.10}
\arraycolsep=1.5pt
\begin{array}[b]{rl}
\D0&\leq\left|(cQ)^\theta(x,t)-(cQ_{\infty})^{\theta}(x)\right|^k\\[5mm]
&\leq\D
k\int_0^x|(cQ)^\theta-(cQ_{\infty})^\theta|^{k-1}|((cQ)^\theta-(cQ_\infty)^\theta)_x|dy\\[5mm]
&\leq\D
k\left(\int^1_0(cQ_{\infty})^{-\eta}|(cQ)^{\theta}-(cQ_{\infty})^{\theta}|^{2k-2}dx\right)^{\frac{1}{2}}
\left(\int^1_0(cQ_{\infty})^{\eta}((cQ)^{\theta}-(cQ_{\infty})^{\theta})^2_xdx\right)^{\frac{1}{2}}\\[5mm]
&=\D
k\left(\int^1_0(cQ_{\infty})^{-\eta}|(cQ)^{\theta}-(cQ_{\infty})^{\theta}|^{\frac{\eta}{\theta}}
|(cQ)^{\theta}-(cQ_{\infty})^{\theta}|^{2k-2-\frac{\eta}{\theta}}dx\right)^{\frac{1}{2}}\\[5mm]
&\D\ \ \ \ \ \  \left(\int^1_0(cQ_{\infty})^{\eta}((cQ)^{\theta}-(cQ_{\infty})^{\theta})^2_xdx\right)^{\frac{1}{2}}\\[5mm]
&\leq\D C\left(\int^1_0
|(cQ)^{\theta}-(cQ_{\infty})^{\theta}|^{2k-2-\frac{\eta}{\theta}}dx\right)^{\frac{1}{2}}
\left(\int^1_0(cQ_{\infty})^{\eta}((cQ)^{\theta}-(cQ_{\infty})^{\theta})^2_xdx\right)^{\frac{1}{2}},
\end{array}
\end{eqnarray}
where $\eta=(2+\alpha)\gamma-\theta-1$.
 Note that
\begin{eqnarray}
|(cQ)^{\theta}-(cQ_{\infty})^{\theta}|\leq
C|cQ-cQ_{\infty}|^{\min\{\theta,1\}}=C|cQ-cQ_{\infty}|^{\theta}.
\end{eqnarray}
Now letting $2k-2-\frac{\eta}{\theta}>0$, we deduce from
\eqref{CQLq} and \eqref{4.10} that
\begin{eqnarray}
\lim_{t\rightarrow\infty}\left|(cQ)^\theta(x,t)-(cQ_{\infty})^{\theta}(x)\right|^k=0,
\end{eqnarray}
uniformly in $x\in[0,1]$, and \eqref{cQlimit} follows. This
completes the proof of Theorem 2.3.

\section{Stabilization rate estimates}
\setcounter{equation}{0} In this section, using the method in
\cite{Fang-Zhang4}, we will give the stabilization rate estimates of
the weak solution $((cQ)(x,t), u(x,t))$ under the condition
$\theta\in(0,\gamma-1)\cap(0,\frac{\gamma}{2}]\cap(0,1-\alpha\gamma]$.
Compared with the corresponding results in
\cite{Fang-Zhang4,Zi-Zhu}, the only difference is that we must deal
with the frictional force $-h(Q)u|u|$.

\bigbreak \noindent\textbf{Lemma 5.1.} Assume the conditions in
Theorem 2.3 hold, we have
\begin{eqnarray}\label{5.1}
\int^1_0\left(u^2+x^{1-\frac{3}{\gamma}}(cQ-cQ_\infty)^2\right)dx\leq\frac{C}{1+t},\
\  \forall \ t\geq 0,
\end{eqnarray}
and
\begin{eqnarray}\label{5.2}
\int^\infty_0\int^1_0(1+t)c^{\theta}Q^{1+\theta}u_x^2dxdt\leq C.
\end{eqnarray}

\bigbreak\noindent{\it Proof.} Multiplying $(\ref{E2})_3$ by
$(1+t)u$, integrating the result equation over $[0,1]\times[0,t]$,
integrating by parts, we get
\begin{eqnarray}
\arraycolsep=1.5pt
\begin{array}[b]{rl}
&\displaystyle
(1+t)\int^1_0\left(\frac{1}{2}u^2+\int^Q_{Q_{\infty}}\frac{c^{\gamma}(h^\gamma-Q^\gamma_\infty)}{\rho_lh^2}dh\right)dx
+\int_0^t\int_0^1(1+s)c^{\theta}Q^{1+\theta}u_x^2dxds\\[5mm]
&\displaystyle +\int^t_0\int^1_0(1+s)h(Q)u^2|u|dxds\\[5mm]
=&\displaystyle
\int^1_0\left(\frac{1}{2}u^2_0+\int^{Q_0}_{Q_{\infty}}\frac{c^{\gamma}(h^\gamma-Q^\gamma_\infty)}{\rho_lh^2}dh\right)dx
+\int^t_0\int^1_0\left(\frac{1}{2}u^2+\int^Q_{Q_{\infty}}\frac{c^{\gamma}(h^\gamma-Q^\gamma_\infty)}{\rho_lh^2}dh\right)dxds \\[5mm]
\leq&\displaystyle
C+C\int^t_0\int^1_0x^{\frac{\beta-1}{\gamma}+\alpha}((cQ)^{\gamma}-(cQ_{\infty})^{\gamma})((cQ_{\infty})^{-\beta}-(cQ)^{-\beta})dxds,
\end{array}
\end{eqnarray}
where we have used \eqref{basic energy}, Lemma 3.4 and the following
fact
$$
\arraycolsep=1.5pt
\begin{array}[b]{rl}
&\displaystyle
\int^t_0\int^1_0\left(\int^Q_{Q_{\infty}}\frac{c^{\gamma}(h^\gamma-Q^\gamma_\infty)}{\rho_lh^2}dh\right)dxds\\[5mm]
\leq&\displaystyle C \int^t_0\int^1_0Q_{\infty}^{-2}|Q-Q_{\infty}||(cQ)^{\gamma}-(cQ_{\infty})^{\gamma}|dxds\\[5mm]
\leq&\displaystyle C \int^t_0\int^1_0c(cQ_{\infty})^{-2}|cQ-cQ_{\infty}||(cQ)^{\gamma}-(cQ_{\infty})^{\gamma}|dxds\\[5mm]
\leq&\displaystyle
C\int^t_0\int^1_0x^{\frac{\beta-1}{\gamma}+\alpha}((cQ)^{\gamma}-(cQ_{\infty})^{\gamma})((cQ_{\infty})^{-\beta}-(cQ)^{-\beta})dxds.
\end{array}
$$
Choose $\beta=1-\alpha\gamma$, then \eqref{5.1} and \eqref{5.2} are
immediately obtained from Lemma 3.8. This completes the proof of
Lemma $5.1$.

\bigbreak \noindent\textbf{Corollary 5.2.} Assume the conditions in
Theorem 2.3 hold, we have
\begin{eqnarray}
\int_0^t \int_0^1(1+s)|u|^idxds\leq C, \ \ \ \  \ \ i=2,3,4.
\end{eqnarray}
\bigbreak\noindent{\it Proof.} The result can be easily obtained by
$(\ref{u2esti})$, $(\ref{u3})$, $(\ref{u4})$ and Lemma 5.1 and the
details are omitted.
\bigbreak \noindent\textbf{Corollary 5.3.}
Assume the conditions in Theorem 2.3 hold, we have
\begin{eqnarray}
\int^t_0\int^1_0(1+s)x^{-\frac{\theta+1}{\gamma}+\alpha}((cQ)^\gamma-(cQ_\infty)^\gamma)^2dxds\leq
C.
\end{eqnarray}
\bigbreak\noindent{\it Proof.} Due to $(\ref{E2})_3$, we have
$$
\arraycolsep=1.5pt
\begin{array}[b]{rl}
&\displaystyle
\int^t_0\int^1_0(1+s)x^{-\frac{\theta+1}{\gamma}+\alpha}((cQ)^\gamma-(cQ_\infty)^\gamma)^2dxds\\[5mm]
=&\displaystyle
\int^t_0\int^1_0(1+s)x^{-\frac{\theta+1}{\gamma}+\alpha}((cQ)^\gamma-(cQ_\infty)^\gamma)\left(-\int^x_0u_tdy\right)dxds\\[5mm]
&+\displaystyle
\int^t_0\int^1_0(1+s)x^{-\frac{\theta+1}{\gamma}+\alpha}((cQ)^\gamma-(cQ_\infty)^\gamma)\left(-\int^x_0h(Q)u|u|dy\right)dxds\\[5mm]
&+\displaystyle
\int^t_0\int^1_0(1+s)x^{-\frac{\theta+1}{\gamma}+\alpha}((cQ)^\gamma-(cQ_\infty)^\gamma)(c^{\theta}Q^{1+\theta}u_x)dxds\\[5mm]
=& J_{1}+J_{2}+J_{3}.
\end{array}
$$
Next, we estimate terms on the right hand side of the above equality
as follows:
\begin{equation}\label{ut}
 \arraycolsep=1.5pt
\begin{array}[b]{rl}
J_{1}=&\displaystyle\int^1_0 x^{-\frac{\theta+1}{\gamma}+\alpha}((cQ_0)^\gamma-(cQ_\infty)^\gamma)\left(\int^x_0u_0dy\right)dx\\[5mm]
&\displaystyle -(1+t)\int^1_0
x^{-\frac{\theta+1}{\gamma}+\alpha}((cQ)^\gamma-(cQ_\infty)^\gamma)\left(\int^x_0udy\right)dx\\[5mm]
&+\displaystyle \int^t_0\int^1_0
x^{-\frac{\theta+1}{\gamma}+\alpha}((cQ)^\gamma-(cQ_\infty)^\gamma)\left(\int^x_0udy\right)dxds\\[5mm]
&-\displaystyle\int^t_0\int^1_0(1+s)x^{-\frac{\theta+1}{\gamma}+\alpha}\gamma\rho_lc^{\gamma}Q^{\gamma+1}u_x\left(\int^x_0udy\right)dxds\\[5mm]
\leq &\displaystyle \int^1_0
x^{\frac{1}{2}-\frac{\theta+1}{\gamma}+1+\alpha}\left(\int^1_0u_0^2dy\right)^{\frac{1}{2}}dx
+(1+t)\int^1_0
x^{\frac{1}{2}-\frac{\theta+1}{\gamma}+1-\frac{2}{\gamma}+\alpha}(cQ-cQ_\infty)^2\left(\int^1_0u^2dy\right)^{\frac{1}{2}}dx\\[5mm]
&\displaystyle+C\int^t_0\int^1_0\left(\int_{0}^xu^2dy\right)dxds+\int^t_0\int^1_0x^{1-2\frac{1+\theta}{\gamma}+2\alpha}((cQ)^\gamma-(cQ_\infty)^\gamma)^2dxds
\\[5mm]
&\displaystyle
+C\int^t_0\int^1_0(1+s)c^{\theta}Q^{1+\theta}u_x^2dxds
+C\int^t_0\int^1_0(1+s)x^{1-2\frac{\theta+1}{\gamma}+2\alpha}c^{2\gamma-\theta}Q^{2\gamma-\theta+1}\left(\int_0^1u^2dy\right)dxds\\[5mm]
\leq & C+C\displaystyle (1+t)\int^1_0
x^{1-\frac{3}{\gamma}}(cQ-cQ_\infty)^2dx+C\int^t_0\int^1_0u^2dxds\\[5mm]
&\displaystyle+\int^t_0\int^1_0x^{2-2\frac{1+\theta}{\gamma}+2\alpha-\frac{\beta}{\gamma}}((cQ)^\gamma-(cQ_\infty)^\gamma)
((cQ)^\beta-(cQ_\infty)^\beta)dxds\\[5mm]
&+\displaystyle
C\int^t_0\int^1_0(1+s)x^{3-\frac{3\theta+1}{\gamma}+\alpha}\left(\int_0^1u^2dy\right)dxds\\[5mm]
\leq & C,
\end{array}
\end{equation}
where we have used \eqref{basic energy}, Lemma 3.4, Lemma 5.1,
Corollary 5.2 and $\beta=2(\gamma-1-\theta+\alpha\gamma)>0$ in Lemma
3.8.
\begin{equation}\label{h(Q)uu}
\arraycolsep=1.5pt
\begin{array}[b]{rl}
J_{2} \leq &\displaystyle
\frac{1}{4}\int^t_0\int^1_0(1+s)x^{-\frac{\theta+1}{\gamma}+\alpha}((cQ)^\gamma-(cQ_\infty)^\gamma)^2dxds\\[5mm]
&\displaystyle
+C\int^t_0\int^1_0(1+s)x^{-\frac{\theta+1}{\gamma}+\alpha}\left(\int^x_0h(Q)u|u|dy\right)^2dxds\\[5mm]
\leq &\displaystyle
\frac{1}{4}\int^t_0\int^1_0(1+s)x^{-\frac{\theta+1}{\gamma}+\alpha}((cQ)^\gamma-(cQ_\infty)^\gamma)^2dxds
+C\int^t_0\int^1_0x^{1-\frac{\theta+1}{\gamma}+\alpha}(1+s)\left(\int_{0}^{x}u^4dy\right)dxds\\[5mm]
\leq &\displaystyle
\frac{1}{4}\int^t_0\int^1_0(1+s)x^{-\frac{\theta+1}{\gamma}+\alpha}((cQ)^\gamma-(cQ_\infty)^\gamma)^2dxds
+C\int^t_0\int^1_0(1+s)u^4dxds\\[5mm]
\leq &\displaystyle
\frac{1}{4}\int^t_0\int^1_0(1+s)x^{-\frac{\theta+1}{\gamma}+\alpha}((cQ)^\gamma-(cQ_\infty)^\gamma)^2dxds+C,
\end{array}
\end{equation}
where we have used Young inequality and the fact $\theta<\gamma-1 $
in Corollary 5.2.

By using Young inequality and Lemma 5.1, we have
\begin{eqnarray}\label{cQux}
\arraycolsep=1.5pt
\begin{array}[b]{rl}
J_{3}\leq &\displaystyle
\frac{1}{4}\int^t_0\int^1_0(1+s)x^{-\frac{\theta+1}{\gamma}+\alpha}((cQ)^\gamma-(cQ_\infty)^\gamma)^2dxds\\[5mm]
&+\displaystyle C\int^t_0\int^1_0(1+s)x^{-\frac{\theta+1}{\gamma}+\alpha}(c^{\theta}Q^{1+\theta})c^{\theta}Q^{1+\theta}u_x^2dxds\\[5mm]
\leq &\displaystyle
\frac{1}{4}\int^t_0\int^1_0(1+s)x^{-\frac{\theta+1}{\gamma}+\alpha}((cQ)^\gamma-(cQ_\infty)^\gamma)^2dxds+C.
\end{array}
\end{eqnarray}
Then from \eqref{ut}, \eqref{h(Q)uu}, \eqref{cQux}, we complete the
proof of Corollary 5.3.

\bigbreak \noindent\textbf{Lemma 5.4.} Assume the conditions in
Theorem 2.3 hold, we have
\begin{eqnarray}\label{5.9}
\int^1_0x^{3-\frac{2\theta+1}{\gamma}+\alpha}\left(((cQ)^\theta)_x-((cQ_\infty)^\theta)_x\right)^2dx\leq
\frac{C}{1+t},\ \  \forall \ t\geq0,
\end{eqnarray}
and
\begin{eqnarray}\label{5.10}
\int^\infty_0\int^1_0(1+t)x^{4-\frac{3\theta+1}{\gamma}+\alpha}\left(((cQ)^\theta)_x-((cQ_\infty)^\theta)_x\right)^2dxdt\leq
C.
\end{eqnarray}
\bigbreak\noindent{\it Proof.} Multiplying (\ref{3.31}) by
$(1+s)x^{3-\frac{2\theta+1}{\gamma}+\alpha}\left(((cQ)^\theta)_x-((cQ_{\infty})^{\theta})_x+\theta\rho_l
u\right)$, and integrating it on $[0, 1]\times[0, t]$, we obtain
\begin{equation}\label{5.11}
\arraycolsep=1.5pt
\begin{array}[b]{rl}
&\displaystyle
\frac{1}{2}(1+t)\int^1_0x^{3-\frac{2\theta+1}{\gamma}+\alpha}\left(((cQ)^\theta)_x-((cQ_{\infty})^{\theta})_x+\theta\rho_l
u\right)^2dx\\[5mm]
&\displaystyle
+\gamma\rho_l\int^t_0\int^1_0(1+s)x^{3-\frac{2\theta+1}{\gamma}+\alpha}(cQ)^{\gamma-\theta}\left(((cQ)^\theta)_x-((cQ_{\infty})^{\theta})_x+\theta\rho_l
u\right)^2dxds\\[5mm]
=&\displaystyle
\frac{1}{2}\int^1_0x^{3-\frac{2\theta+1}{\gamma}+\alpha}\left(((cQ_0)^\theta)_x-((cQ_{\infty})^{\theta})_x+\theta\rho_l
u_0\right)^2dx\\[5mm]
&\displaystyle
+\frac{1}{2}\int^t_0\int^1_0x^{3-\frac{2\theta+1}{\gamma}+\alpha}\left(((cQ)^\theta)_x-((cQ_{\infty})^{\theta})_x+\theta\rho_l
u\right)^2dxds\\[5mm]
&\displaystyle
+\int^t_0\int^1_0(1+s)x^{3-\frac{2\theta+1}{\gamma}+\alpha}\left(((cQ)^\theta)_x-((cQ_{\infty})^{\theta})_x+
\theta\rho_l
u\right)\theta\rho_{l}^2\gamma(cQ)^{\gamma-\theta}udxds\\[5mm]
&\displaystyle+\theta\rho_{l}
g\int^t_0\int^1_0(1+s)x^{3-\frac{2\theta+1}{\gamma}+\alpha}\left(((cQ)^\theta)_x-((cQ_{\infty})^{\theta})_x+
\theta\rho_l
u\right)\left(1-\frac{(cQ)^{\gamma-\theta}}{(cQ_{\infty})^{\gamma-\theta}}\right)dxds\\[5mm]
&\displaystyle+\int^t_0\int^1_0(1+s)x^{3-\frac{2\theta+1}{\gamma}+\alpha}\left(((cQ)^\theta)_x-((cQ_{\infty})^{\theta})_x+\theta\rho_l
u\right)(-\theta\rho_{l} h(Q)u|u|)dxds\\[5mm]
\leq&\displaystyle
C+\frac{1}{2}\int^t_0\int^1_0(1+s)x^{3-\frac{2\theta+1}{\gamma}+\alpha}\gamma\rho_{l}(cQ)^{\gamma-\theta}
\left(((cQ)^\theta)_x-((cQ_{\infty})^{\theta})_x+\theta\rho_{l}
u\right)^2dxds\\[5mm]
&\displaystyle
+C\int^t_0\int^1_0(1+s)x^{3-\frac{2\theta+1}{\gamma}+\alpha}(cQ)^{\gamma-\theta}u^2dxds\\[5mm]
&\displaystyle+C\int^t_0\int^1_0(1+s)x^{3-\frac{2\theta+1}{\gamma}+\alpha}(cQ)^{\theta-\gamma}
\left(1-\frac{(cQ)^{\gamma-\theta}}{(cQ_{\infty})^{\gamma-\theta}}\right)^2dxds\\[5mm]
&+\displaystyle
C\int^t_0\int^1_0(1+s)x^{3-\frac{2\theta+1}{\gamma}+\alpha}(cQ)^{\theta-\gamma}h^2(Q)u^4dxds,
\end{array}
\end{equation}
where we have used Lemma 3.9 and the assumption $(A_{2})$.

 By using Lemma 3.3, Corollaries 5.1 and 5.3, we get
\begin{equation}\label{5.12}
\arraycolsep=1.5pt
\begin{array}[b]{rl}
&\displaystyle
\int^t_0\int^1_0(1+s)x^{3-\frac{2\theta+1}{\gamma}+\alpha}(cQ)^{\theta-\gamma}
\left(1-\frac{(cQ)^{\gamma-\theta}}{(cQ)^{\gamma-\theta}_{\infty}}\right)^2dxds\\[5mm]
\leq&\displaystyle
C\int^t_0\int^1_0(1+s)x^{-\frac{\theta+1}{\gamma}+\alpha}((cQ)^\gamma-(cQ)^\gamma_\infty)^2dxds\leq
C,
\end{array}
\end{equation}
and
\begin{equation}\label{5.13}
\int^t_0\int^1_0(1+s)x^{3-\frac{2\theta+1}{\gamma}+\alpha}(cQ)^{\theta-\gamma}h^2(Q)u^4dxds\leq
C\int^t_0\int^1_0(1+s)u^4dxds\leq C.
\end{equation}
From \eqref{5.11}-\eqref{5.13}, we obtain \eqref{5.9} and
\eqref{5.10} and the proof of Lemma 5.4 is complete.

\bigbreak Following the same ideas in Proposition 6.3 of
\cite{Fang-Zhang4}, we can estimate the stabilization rate of
$u(x,t)$ in the sense of $L^\infty$ norm under the conditions of
Theorem 2.4.

\bigbreak\noindent\textbf{Lemma 5.5.} Assume the conditions in
Theorem 2.3 hold, we have
\begin{equation}\label{5.14}
\int^1_0c^{\theta}Q^{1+\theta}u^2_xdx\leq\frac{C}{1+t}, \ \  \forall \
t\geq0,
\end{equation}
and
\begin{equation}
\int^\infty_0\int^1_0(1+t)u^2_tdxdt \leq C.
\end{equation}
\bigbreak\noindent{\it Proof.} Multiplying $(\ref{E2})_3$ by
$(1+s)u_t$, integrating the result equation over $[0,1]\times[0,t]$,
integrating by parts, we get
\begin{equation}\label{5.16}
\arraycolsep=1.5pt
\begin{array}[b]{rl}
&\displaystyle\int_0^t\int_0^1(1+s)u_t^2dxds-\int_0^t\int_0^1(1+s)(cQ)^{\gamma}u_{tx}dxds\\[5mm]
=&\displaystyle-\frac{1}{3}(1+t)\int_0^1h(Q)|u|^3dx+\frac{1}{3}\int_0^1h(Q_0)|u_0|^3dx\\[5mm]
&\displaystyle -\frac{1}{3}\rho_{l}\int_0^t\int_0^1(1+s)h'(Q)Q^2|u|^3u_xdxds+\frac{1}{3}\int_0^t\int_0^1h(Q)|u|^3dxds\\[5mm]
&\displaystyle
-\int_0^t\int_0^1(1+s)(cQ_{\infty})^{\gamma}u_{tx}dxds
-\int_0^t\int_0^1(1+s)(c^{\theta}Q^{1+\theta}u_x)u_{tx}dxds.
\end{array}
\end{equation}
The last term on the right-hand side in \eqref{5.16} can be
estimated as follows:
\begin{equation}\label{5.17}
\arraycolsep=1.5pt
\begin{array}[b]{rl}
&\displaystyle-\int_0^t\int_0^1(1+s)(c^{\theta}Q^{1+\theta}u_x)u_{tx}dxds\\[5mm]
=&\displaystyle
-\frac{1}{2}(1+t)\int_0^1c^{\theta}Q^{1+\theta}u_x^2dx+\frac{1}{2}\int_0^1c^{\theta}Q_0^{1+\theta}u_{0x}^2dx\\[5mm]
&+\displaystyle
\frac{1}{2}\int_0^t\int_0^1c^{\theta}Q^{1+\theta}u_x^2dxds-\frac{1}{2}(1+\theta)\rho_{l}\int_0^t\int_0^1(1+s)c^{\theta}Q^{2+\theta}u_x^3dxds.
\end{array}
\end{equation}
Substituting \eqref{5.17} into \eqref{5.16}, we get
\begin{equation}\label{5.18}
\arraycolsep=1.5pt
\begin{array}[b]{rl}
&\displaystyle
\int_0^t\int_0^1(1+s)u_t^2dxds+\frac{1}{2}(1+t)\int_0^1c^{\theta}Q^{1+\theta}u_x^2dx+\frac{1}{3}(1+t)\int_0^1h(Q)|u|^3dx\\[5mm]
=&\displaystyle
\frac{1}{3}\int_0^1h(Q_0)|u_0|^3dx-\frac{1}{3}\rho_{l}\int_0^t\int_0^1(1+s)h'(Q)Q^2|u|^3u_xdxds
+\frac{1}{3}\int_0^t\int_0^1h(Q)|u|^3dxds\\[5mm]
&\displaystyle+\frac{1}{2}\int_0^1c^{\theta}Q_0^{1+\theta}u_{0x}^2dx
+\frac{1}{2}\int_0^t\int_0^1c^{\theta}Q^{1+\theta}u_x^2dxds-\frac{1}{2}(1+\theta)\rho_{l}\int_0^t\int_0^1(1+s)c^{\theta}Q^{2+\theta}u_x^3dxds\\[5mm]
&\D
+\int_0^t\int_0^1(1+s)((cQ)^{\gamma}-(cQ_{\infty})^{\gamma})u_{tx}dxds\\[5mm]
=& \displaystyle\sum_{i=1}^{i=7}K_{i}.
\end{array}
\end{equation}
It is easy to see $K_{1}+K_{3}+K_{4}+K_{5}\leq C$ from the
assumption $(A_{2})$, Lemma 3.2 and Lemma 3.5.
\begin{equation}\label{5.22}
\arraycolsep=1.5pt
\begin{array}{rl}
\displaystyle
K_{2}\leq &\displaystyle C\int_0^t\int_0^1(1+s)\frac{Q^3}{(1+Q)^3}|u|^3|u_x|dxds\\[5mm]
\leq & \displaystyle \|Qu_{x}\|_{L^\infty}
\int_0^t\int_0^1(1+s)|u|^3dxds\leq C,
\end{array}
\end{equation}
where we have used $\|Qu_{x}\|_{L^\infty}\leq C$ and Corollary 5.2.

By using $\|Qu_{x}\|_{L^\infty}\leq C$ and Lemma 5.1, we have
\begin{equation}\label{5.23}
K_{6} \leq
C\|Qu_{x}\|_{L^\infty}\int_0^t\int_0^1(1+s)c^{\theta}Q^{1+\theta}u_x^2dxds\leq
C.
\end{equation}
In order to complete the proof of Lemma 5.5, it suffices to
 estimate $K_{7}$ of the right-hand side in \eqref{5.18}.
\begin{equation}\label{J7}
\arraycolsep=1.5pt
\begin{array}[b]{rl}
K_{7}=&\displaystyle(1+t)\int_0^1((cQ)^{\gamma}-(cQ_{\infty})^{\gamma})u_xdx-\int_0^1((cQ_0)^{\gamma}-(cQ_{\infty})^{\gamma})u_{0x}dx\\[5mm]
&\displaystyle
+\gamma\rho_l\int_0^t\int_0^1(1+s)c^{\gamma}Q^{1+\gamma}u_x^2dxds-\int_0^t\int_0^1((cQ)^{\gamma}-(cQ_{\infty})^{\gamma})u_xdxds\\[5mm]
\leq&\displaystyle
\frac{1}{4}(1+t)\int_0^1c^{\theta}Q^{1+\theta}u_x^2dx+
C(1+t)\int_0^1((cQ)^{\gamma}-(cQ_{\infty})^{\gamma})^2c^{-\theta}Q^{-1-\theta}dx\\[5mm]
&+\displaystyle\int_0^1((cQ_0)^{\gamma}-(cQ_{\infty})^{\gamma})^2dx+\int_0^1u_{0x}^{2}dx\\[5mm]
&\displaystyle
+\gamma\rho_{l}\int_0^t\int_0^1(cQ)^{\gamma-\theta}(1+s)c^{\theta}Q^{1+\theta}u_x^2dxds
+\int_0^t\int_0^1((cQ)^{\gamma}-(cQ_{\infty})^{\gamma})^2c^{-\theta}Q^{-1-\theta}dxds\\[5mm]
&\displaystyle +\int_0^t\int_0^1c^{\theta}Q^{1+\theta}u_x^2dxds\\[5mm]
\leq & \displaystyle C
+C(1+t)\int_0^1((cQ)^{\gamma}-(cQ_{\infty})^{\gamma})^2c^{-\theta}Q^{-1-\theta}dx\\[5mm]
&+\displaystyle\int_0^t\int_0^1((cQ)^{\gamma}-(cQ_{\infty})^{\gamma})^2c^{-\theta}Q^{-1-\theta}dxds.
\end{array}
\end{equation}
Here, we have used Lemma 5.1, the assumption $(A_{2})$ and Lemma
3.2. The rest two terms on the right-hand side of \eqref{J7} can be
estimated as follows:
\begin{equation}\label{5.20}
\arraycolsep=1.5pt
\begin{array}[b]{rl}
\displaystyle
C(1+t)\int_0^1((cQ)^{\gamma}-(cQ_{\infty})^{\gamma})^2c^{-\theta}Q^{-1-\theta}dx
\leq &\displaystyle
C(1+t)\int_0^1(cQ-cQ_{\infty})^2x^{2-\frac{3+\theta}{\gamma}}dx\\[5mm]
\leq &\displaystyle
C(1+t)\int_0^1(cQ-cQ_{\infty})^2x^{1-\frac{3}{\gamma}}dx\\[5mm]
\leq &\displaystyle C,
\end{array}
\end{equation}
provided $\theta\leq \gamma-1$. And if we take
$\beta=\gamma-\theta-1$ in Lemma 3.9, we have
\begin{equation}\label{5.21}
\arraycolsep=1.5pt
\begin{array}[b]{rl}
&\displaystyle\int_0^t\int_0^1((cQ)^{\gamma}-(cQ_{\infty})^{\gamma})^2c^{-\theta}Q^{-1-\theta}dxds\\[5mm]
=&\displaystyle\int_0^t\int_0^1\frac{(cQ)^{\gamma}-(cQ_{\infty})^{\gamma}}{(cQ)^{\beta}-(cQ_{\infty})^{\beta}}c^{-\theta}Q^{-1-\theta}((cQ)^{\gamma}-(cQ_{\infty})^{\gamma})
((cQ)^{\beta}-(cQ_{\infty})^{\beta})dxds\\[5mm]
\leq & \displaystyle
C\int_0^t\int_0^1((cQ)^{\gamma}-(cQ_{\infty})^{\gamma})((cQ)^{\beta}-(cQ_{\infty})^{\beta})dxds
\leq C.
\end{array}
\end{equation}
Finally, from \eqref{5.18}-\eqref{5.21} we have
$$
\int_0^t\int_0^1(1+s)u_t^2dxds+(1+t)\int_0^1c^{\theta}Q^{1+\theta}u_x^2dx+(1+t)\int_0^1h(Q)|u|^2dx\leq
C.
$$
This proves Lemma 5.5.

\bigbreak {\noindent\underline{\it Proof of Theorem 2.4}}.

Choosing a positive number $m$ large enough, which is to be
determined later, applying H\"{o}lder inequality, \eqref{2.6} and
\eqref{5.9},  we have by $(cQ)(0,t)=(cQ_{\infty})(0,t)=0$
\begin{equation}\label{5.24}
\arraycolsep=1.5pt
\begin{array}[b]{rl}
\displaystyle0&\leq\left|(cQ)^\theta(x,t)-(cQ_{\infty})^{\theta}(x,t)\right|^m\\[5mm]
&\leq\displaystyle
m\int_0^x|(cQ)^\theta-(cQ_\infty)^\theta|^{m-1}|((cQ)^\theta-(cQ_\infty)^\theta)_x|dy\\[5mm]
&\leq\displaystyle
m\left(\int^1_0x^{-\alpha+\frac{2\theta+1}{\gamma}-3}|(cQ^{\theta}-(cQ_{\infty})^{\theta}|^{2m-2}dx\right)^{\frac{1}{2}}
\left(\int^1_0x^{3-\frac{2\theta+1}{\gamma}+\alpha}((cQ)^{\theta}-(cQ_{\infty})^{\theta})^2_xdx\right)^{\frac{1}{2}}\\[5mm]
&\leq\displaystyle
C(1+t)^{-\frac{1}{2}}\left(\int^1_0x^{-\alpha+\frac{2\theta+1}{\gamma}-3}
|(cQ)^{\theta}-(cQ_{\infty})^{\theta}|^{2m-2}dx\right)^{\frac{1}{2}}\\[5mm]
&\leq\displaystyle
C(1+t)^{-\frac{1}{2}}\left(\int^1_0x^{\frac{-\alpha\gamma+2\theta+2-4\gamma+\theta(2m-2)}{\gamma}}
\left(x^{1-\frac{3}{\gamma}}(cQ-cQ_{\infty})^2\right)dx\right)^{\frac{1}{2}}.
\end{array}
\end{equation}
Taking $m=\frac{4\gamma+\alpha\gamma-2}{2\theta}$, then
$2\theta+2-4\gamma+\theta(2m-2)-\alpha\gamma=0$. It follows from
\eqref{5.24} and \eqref{5.1} that \eqref{2.12} holds.

For the velocity function $u(x,t)$, \eqref{2.13} follows from
(\ref{4.1}) and (\ref{5.14}) directly. This completes the proof of
Theorem 2.4.

\bigbreak

\noindent{\bf Acknowledgements:}\ \  The authors were supported by
the Natural Science Foundation of China $\#$11071093, the PhD
specialized grant of the Ministry of Education of China
$\#$20100144110001, and the Special Fund for Basic Scientific
Research  of Central Colleges $\#$CCNU10C01001.

\bigbreak

\bibliographystyle{plain}

\end{document}